\documentclass[11pt]{article}

\usepackage[symbol]{footmisc}

\usepackage[
  shownumpages,               % comment to hide page numbers
  bgcolor={245,245,250},      % box bg color: RGB or xcolor name (e.g. gray, blue, etc.)
  braincolor={0,138,255},    % accent color: RGB or xcolor name (e.g. gray, blue, etc.)
  linkcolor={0,62,165},  % defaults to braincolor
  citecolor={0,62,165},  % defaults to braincolor
  urlcolor={0,62,165},            % defaults to blue
  citingstyle=authoryear,     % authoryear: [Ivanov et al., 2023] | numbers: [1] | super: ¹
  bibliostyle=plainnat,       % plainnat, abbrvnat, unsrtnat, ...
  bibfile=brainlab          % .bib file name (without .bib)
]{brainlab}

\definecolor{mydarkred}{RGB}{192,47,25}
\definecolor{mydarkgreen}{RGB}{39,130,67}
\newcommand{\green}{\color{mydarkgreen}}
\newcommand{\red}{\color{mydarkred}}

\newcommand{\cmark}{{\green\ding{51}}}
\newcommand{\xmark}{{\red\ding{55}}}

\setbrainmeta{
  title={Shuffling Heuristic in Variational Inequalities: Establishing New Convergence Guarantees},
  authors={
    Daniil Medyakov\textsuperscript{1,2}\footnotemark, Gleb Molodtsov\textsuperscript{1,2}, Grigoriy Evseev\textsuperscript{1}, Egor Petrov\textsuperscript{1}, Aleksandr Beznosikov\textsuperscript{1,2}
  },
  affiliations={
    \textsuperscript{1}Basic Research of Artificial Intelligence Laboratory (BRAIn Lab) \\
    \textsuperscript{2}Federated Learning Problems Laboratory
  },
  abstract={
    Variational inequalities have gained significant attention in machine learning and optimization research. While stochastic methods for solving these problems typically assume independent data sampling, we investigate an alternative approach -- the shuffling heuristic. This strategy involves permuting the dataset before sequential processing, ensuring equal consideration of all data points. Despite its practical utility, theoretical guarantees for shuffling in variational inequalities remain unexplored. We address this gap by providing the first theoretical convergence estimates for shuffling methods in this context. Our analysis establishes rigorous bounds and convergence rates, extending the theoretical framework for this important class of algorithms. We validate our findings through extensive experiments on diverse benchmark variational inequality problems, demonstrating faster convergence of shuffling methods compared to independent sampling approaches.
  },
}

\begin{document}

\begin{mainpart}

\footnotetext{Corresponding author: medyakovd3@gmail.com}

\section{Introduction}

Variational inequalities (VIs) have attracted researchers' attention in various fields for more than half a century \cite{browder1965nonexpansive}. In this work, we investigate the variational inequality problem in the following form:
\begin{equation}\label{eq:vi_setting}
    \text{find}~~ z^* \in Z ~~ \text{such that}~~ \forall z \in Z \hookrightarrow\langle F(z^*), z - z^*\rangle + g(z) - g(z^*) \geqslant 0,
\end{equation}
where \( F \) is a monotone operator and \( g \) is a proper convex lower semicontinuous function, which serves as a regularizer. Variational inequalities serve as a universal tool for addressing specific problems, such as minimization, saddle point problems, fixed point problems, and others \cite{facchinei2003finite,kinderlehrer2000introduction}. We provide examples to offer intuition about VIs. 
\begin{example}[Convex optimization]\label{ex:convopt}
    Consider the following convex regularized optimization problem:
    \begin{equation}\label{eq:convregprob}
        \underset{z\in\mathbb R^d}{\min} \left[f(z) + g(z)\right].
    \end{equation}
    In this example, $f$ is a smooth data representative term, and $g$ is probably a non-smooth regularizer. We define $F(z) = \nabla f(z)$. Then $z^*\in \text{dom}~ g$ is the solution of \eqref{eq:vi_setting} if and only if $z^*\in \text{dom~} g$ is the solution of \eqref{eq:convregprob}. In this way, the problem \eqref{eq:convregprob} can be considered as a variational inequality.
\end{example}

\begin{example}[Convex-concave saddles]\label{ex:convconcsaddle}
    Consider the following convex-concave saddle point problem:
    \begin{equation}\label{eq:convconcsaddle}
        \underset{x\in \mathbb R^{d_x}}{\min}\underset{y\in \mathbb R^{d_y}}{\max} \left[f(x,y) + g_1(x) - g_2(y)\right].
    \end{equation}
    There, $f$ has the same interpretation as in Example \ref{ex:convopt}, and $g_1$, $g_2$ can also be perceived as regularizers. We define $F(z) = F(x, y) = \left[\nabla_x f(x, y), -\nabla_y f(x, y)\right]$. Then $z^*\in \text{dom}~ g_1 \times \text{dom}~ g_2$ is the solution of \eqref{eq:vi_setting} if and only if $z^*\in \text{dom}~ g_1 \times \text{dom}~ g_2$ is the solution of \eqref{eq:convconcsaddle}. In this way, the problem \eqref{eq:convconcsaddle} can be considered as a variational inequality.
\end{example}

There are multiple practical reasons to focus on this formulation. Firstly, for numerous non-smooth problems, solutions are obtained more efficiently when formulated as saddle point problems \cite{nesterov2005smooth,nemirovski2004prox,chambolle2011first,esser2010general}. Secondly, recent studies have established new connections between VIs and reinforcement learning \cite{omidshafiei2017deep,jin2020efficiently}, adversarial training \cite{madry2017towards}, and generative adversarial networks (GANs) \cite{goodfellow2014generative}. In particular, consideration of monotone and strongly monotone inequalities provides useful methods and recommendations for the GAN community \cite{daskalakis2017training,gidel2018variational,mertikopoulos2018optimistic,chavdarova2019reducing,liang2019interaction,peng2020training}. VIs also have extensive applications in various classic problems, including discriminative clustering \cite{xu2004maximum}, matrix factorization \cite{bach2008convex}, image denoising \cite{esser2010general,chambolle2011first}, robust optimization \cite{ben2009robust}, economics, game theory \cite{von1953theory}, and optimal control \cite{facchinei2003finite}. 

Solving the problem \eqref{eq:vi_setting} requires specialized methods, as traditional optimization techniques, e.g. the \textsc{Gradient Descent} method, often fall short when applied to VIs and saddle point problems \cite{harker1990finite}. These classic methods not only struggle with efficiency but also provide weak theoretical convergence guarantees in the VI context \cite{beznosikov2023smooth}. Among the various approaches developed for VIs, the \textsc{Extragradient} method \cite{korpelevich1976extragradient,mokhtari2020unified} stands out as one of the most fundamental and effective techniques.

While variational inequalities provide a powerful framework for addressing a wide range of problems, recent trends in machine learning and data science present new challenges. The exponential growth in dataset sizes and the increasing complexity of models create a pressing need for more efficient computational approaches \cite{bottou2010large,dean2012large,medyakov2023optimal}. To address these challenges in the context of VIs, we reformulate the problem by considering the operator $F$ as the finite sum of operators $F_i$:
\begin{equation}\label{eq:finite-sum}
    F(z) = \frac{1}{n}\sum\nolimits_{i=1}^n F_i(z),
\end{equation}
where each $F_i$ corresponds to an individual data point. This decomposition enables us to address large-scale problems more effectively.

In this paper, we explore stochastic algorithms that are particularly suitable for practical extensive applications. As mentioned previously, the number of operators $n$ is typically large, making the computation of the full operator value at each iteration computationally expensive. Instead, stochastic algorithms randomly select $F_i$ at each iteration. The stochastic version of the \textsc{Extragradient} method \cite{juditsky2011solving} selects random independent indices $i_t, j_t$ at iteration $t$ and performs the following updates:
\begin{equation*}
    \begin{aligned}
        z^{t+\frac{1}{2}} &= z^t - \gamma F_{i_t} (z^t),\\
        z^{t+1} &= z^t - \gamma F_{j_t} (z^{t + \frac{1}{2}}).
    \end{aligned}
\end{equation*}
\vspace{-4mm}

Just as deterministic \textsc{Extragradient} is a modification of the classic gradient method with an additional step, the stochastic \textsc{Extragradient} represents the same modification of \textsc{SGD} \cite{robbins1951stochastic}. Although this method performs well on variational inequalities, it encounters a significant issue with its properties and performance that has been thoroughly studied: the variance of its inherent stochastic estimators of operators remains high throughout the learning process. Hence, \textsc{Extragradient} with a constant learning rate converges linearly only to a neighborhood of the optimal solution, the size of which is proportional to the step size and variance \cite{juditsky2011solving}. This problem is also characteristic of classic \textsc{SGD} \cite{bottou2009curiously,moulines2011non,gower2020variance}. 

To address this limitation, the variance reduction (VR) technique has been developed for a classic finite-sum minimization task \cite{johnson2013accelerating}. The method involves the following steps: at the $t$-th iteration, an index $i_t$ is selected along with a reference point $\omega^t$, which is updated once per epoch or selected probabilistically (as in loopless versions, e.g., \cite{kovalev2020don}). Considering the convex optimization problem (see Example \ref{ex:convopt}), we can formally write the stochastic reduced gradient at the point $z^{t+\frac{1}{2}}$ as
\begin{equation*}
    \nabla\hat{f}_{i_t} (z^{t+\frac{1}{2}}) = \nabla f_{i_t} (z^{t+\frac{1}{2}}) - \nabla f_{i_t} (\omega^t) + \nabla f (\omega^t).
\end{equation*} 
\vspace{-4mm}

The objective of the variance reduction mechanisms is to overcome the limitations of naive gradient estimators. These mechanisms employ an iterative process to construct and apply a gradient estimator with progressively reduced variance. This approach permits the use of larger learning rates, thereby accelerating the training process. 

Besides, along with the aforementioned \textsc{SVRG}, there are popular methods for solving the classic finite-sum problem based on this technique, such as \textsc{SAG} \cite{roux2012stochastic}, \textsc{SAGA} \cite{defazio2014saga}, \textsc{Finito} \cite{defazio2014finito}, \textsc{SARAH} \cite{nguyen2017sarah,hu2019efficient}, and \textsc{SPIDER} \cite{fang2018spider}. The variance reduction mechanism is utilized not only in methods that address the minimization problem but also in methods for the problem \eqref{eq:vi_setting}. Examples include the VR versions of \textsc{Extragradient}, \textsc{Mirror-prox} \cite{alacaoglu2022stochastic}, \textsc{gradient method} \cite{palaniappan2016stochastic}, and \textsc{forward-reflected-backward (FoRB)} \cite{alacaoglu2021forward}.

In addition to stochastic methods, various heuristics for selecting the \(i_t\)-th index at each iteration of algorithms exist. Careful analysis of these strategies may lead to the development of more robust and efficient algorithms. In this paper, we explore the shuffling heuristic \cite{mishchenko2020random,safran2020good,koloskova2023convergence,malinovsky2023random}. Unlike the random and independent selection of the index \(i_t\) at each iteration, which is common in classic stochastic methods, this heuristic adopts a more practical approach. Specifically, it involves permuting the sequence of indices \(\{1, \ldots, n\}\), where $n$ is the number of data samples \eqref{eq:finite-sum}, and then selecting the index corresponding to the iteration number during the algorithm's execution. This approach ensures that during one epoch of training, we take a step for each operator, and only once. There are several shuffling techniques available. Among the most popular are Random Reshuffling (RR) \cite{gurbuzbalaban2021random,haochen2019random,nagaraj2019sgd}, where data is shuffled before each epoch; Shuffle Once (SO) \cite{safran2020good,rajput2020closing}, where shuffling occurs once before the start of training; and Cyclic permutation \cite{mangasarian1993serial,bertsekas2000gradient,nedic2001incremental,li2019incremental}, where data is accessed deterministically in a cyclic order.

\textbf{Related Works.}
There are many methods available to solve the problem of variational inequalities. As mentioned above, the standard deterministic choice for solving the problem \eqref{eq:vi_setting} is \textsc{Extragradient} \cite{korpelevich1976extragradient}. This method addresses variational inequalities in the Euclidean setup. Later, \textsc{Mirror-prox} \cite{nemirovski2004prox}, which exploits the Bregman divergence, was proposed. This approach accounts for generalized geometry that may be non-Euclidean. Additionally, there is a set of deterministic methods for solving variational inequalities: \textsc{forward-backward-forward (FBF)} \cite{tseng2000modified}, \textsc{Dual extrapolation} \cite{nesterov2007dual}, \textsc{reflected gradient} \cite{malitsky2015projected}, \textsc{forward-reflected-backward (FoRB)} \cite{malitsky2020forward}. 

For the first time, the stochastic version of algorithms for solving VIs was proposed in the work \cite{juditsky2011solving}. Later, to reduce the variance inherent in these stochastic methods, researchers adopted variance reduction techniques. The initial works in this field are \cite{palaniappan2016stochastic,chavdarova2019reducing}. In particular, \cite{palaniappan2016stochastic} studied the stochastic \textsc{gradient method} with variance reduction. The method was based on \textsc{SVRG} \cite{johnson2013accelerating} and incorporated Catalyst envelope acceleration. The combination of \textsc{Extragradient} and \textsc{SVRG} was considered in \cite{chavdarova2019reducing}. They achieved a worse convergence rate compared to \cite{palaniappan2016stochastic} and only in the strongly monotone setting. Consequently, a notable paper in which the authors considered monotone operators was presented \cite{carmon2019variance}. This work also falls under the Bregman setup but requires additional assumptions on the operator $F$
and considers the matrix games setup. The current state-of-the-art in this area is the article \cite{alacaoglu2022stochastic}, which improved the convergence estimates of all previous studies. This work addressed various scenarios, including generally monotone and strongly monotone operators, as well as the Bregman and Euclidean setups.
Convergence results from the papers highlighted above are summarized in Table \ref{tab:methodscomp}.

%\begin{figure}[htbp]
%\begin{minipage}[t]{1\linewidth}
\begin{table} 
    \centering
    \small
\resizebox{\textwidth}{!}{%
\begin{threeparttable}
\caption{\label{tab:methodscomp} Comparison of the convergence results for the methods for solving VI.}
\renewcommand{\arraystretch}{2.5}
\begin{tabular}{|c|c|c|c|c|}
\hline
\textbf{Algorithm} & \textbf{Sampling} & \textbf{VR?} &
\begin{tabular}{@{}l@{}}
\vspace{-3mm}\textbf{Strongly} \\
\vspace{-3mm}\textbf{Monotone} \\  \textbf{Complexity}
\end{tabular}
 & \begin{tabular}{@{}l@{}}
\vspace{-3mm} \textbf{Monotone} \\  \textbf{Complexity}
\end{tabular} \\ \hline
Extragradient \cite{korpelevich1976extragradient,mokhtari2020unified} & Deterministic &  \xmark & $\mathcal{\widetilde{O}}\left(\frac{nL}{\mu}\right)$ &  $\mathcal{O}\left(\frac{nL}{\varepsilon}\right)$ \\ \hline
Mirror-prox \cite{nemirovski2004prox} & Deterministic & \xmark & $\backslash$ &  $\mathcal{O}\left(\frac{nL}{\varepsilon}\right)$ \\ \hline
FBF \cite{tseng2000modified} & Deterministic & \xmark & $\backslash$ &  $\mathcal{O}\left(\frac{nL}{\varepsilon}\right)$ \\ \hline
FoRB \cite{malitsky2020forward} & Deterministic & \xmark & $\backslash$ &  $\mathcal{O}\left(\frac{nL}{\varepsilon}\right)$ \\ \hline
Mirror-prox \cite{juditsky2011solving} & Independent & \xmark & $\backslash$ &  $\mathcal{O}\left(\frac{L}{\varepsilon} + \frac{1}{\varepsilon^2}\right)$ \\ \hline
Extragradient \cite{beznosikov2020distributed} & Independent & \xmark & $\mathcal{\widetilde{O}}\left(\frac{L}{\mu} + \frac{1}{\mu^2\varepsilon}\right)$ &  $\mathcal{O}\left(\frac{L}{\varepsilon} + \frac{1}{\varepsilon^2}\right)$ \\ \hline
REG \cite{mishchenko2020revisiting} & Independent & \xmark & $\mathcal{\widetilde{O}}\left(\frac{L}{\mu} + \frac{1}{\mu^2\varepsilon}\right)$&  $\mathcal{O}\left(\frac{L}{\varepsilon} + \frac{1}{\varepsilon^2}\right)$  \\ \hline
Extragradient \cite{carmon2019variance} & Independent & \cmark & $\backslash$  &  $\mathcal{\widetilde{O}}\left(n + \frac{\sqrt{n}\overline{L}}{\varepsilon}\right)$ \\ \hline
Mirror-prox \cite{carmon2019variance} & Independent & \cmark & $\backslash$ &  $\mathcal{\widetilde{O}}\left(n + \frac{\sqrt{n}\overline{L}}{\varepsilon}\right)$ \\ \hline
FBF \cite{palaniappan2016stochastic} & Independent & \cmark & $\mathcal{\widetilde{O}}\left(n + \frac{\sqrt{n}\overline{L}}{\mu}\right)$ &  $\mathcal{\widetilde{O}}\left(n + \frac{\sqrt{n}\overline{L}}{\varepsilon}\right)^{\red{(1)}}$ \\ \hline
Extragradient \cite{chavdarova2019reducing} & Independent & \cmark & $\mathcal{\widetilde{O}}\left(n + \frac{\overline{L}^2}{\mu^2}\right)$ &  $\mathcal{\widetilde{O}}\left(n + \frac{\overline{L}^2}{\varepsilon^2}\right)^{\red{(1)}}$ \\ \hline
FoRB \cite{alacaoglu2021forward} & Independent & \cmark & $\backslash$&  $\mathcal{O}\left(n + \frac{n\overline{L}}{\varepsilon}\right)$ \\ \hline
Extragradient \cite{alacaoglu2022stochastic} & Independent & \cmark & $\mathcal{\widetilde{O}}\left(n + \frac{\sqrt{n}\overline{L}}{\mu}\right)$ &  $\mathcal{O}\left(n + \frac{\sqrt{n}\overline{L}}{\varepsilon}\right)$ \\ \hline
Mirror-prox \cite{alacaoglu2022stochastic} & Independent & \cmark & $\backslash$ &  $\mathcal{O}\left(n + \frac{\sqrt{n}\overline{L}}{\varepsilon}\right)$ \\ \hline
\rowcolor{yellow} Extragradient (this paper) & RR / SO & \xmark & $\mathcal{\widetilde{O}}\left(n + \frac{L}{\mu} + \frac{n^2}{\mu^2\varepsilon}\right)$ & $\mathcal{\widetilde{O}}\left(n + \frac{L}{\varepsilon} + \frac{n^2}{\varepsilon^3}\right)^{\red{(1)}}$ \\ \hline
\rowcolor{yellow} Extragradient (this paper) & RR / SO & \cmark & $\mathcal{\widetilde{O}}\left(n\frac{L^2}{\mu^2}\right)$ &  $\mathcal{\widetilde{O}}\left(n\frac{L^2}{\varepsilon^2}\right)^{\red{(1)}}$  \\ \hline 
\end{tabular}%
\begin{tablenotes}
    \item [] {\em Columns:} Sampling = Deterministic, if considered non-stochastic method, Independent, if method uses independent choice of operator's indices, RR / SO if method uses shuffling heuristic, Assumption = assumption on operator $F$, VR? = whether the method uses variance reduction technique.
    \item [] {\em Notation:} $\mu$ = constant of strong monotonicity, $L$ = Lipschitz constant of $F$, $\overline{L}$ = Lipschitz in mean constant, i.e. $\nicefrac{1}{n}\sum\nolimits_{i=1}^n \|F_i(z_1) - F_i(z_2)\| \leqslant L \|z_1 - z_2\| ~\forall z_1, z_2 \in Z$, $n$ =  size of the dataset, $\varepsilon$ = accuracy of the solution.
    \item [] {\red{(1)}}: This result is obtained with regularization trick: $\mu \sim \nicefrac{\varepsilon}{D^2}$.
\end{tablenotes}    
\end{threeparttable}
}
\end{table}
%\end{minipage}
%\end{figure}

\vspace{-1mm}
In all these papers, the estimates were obtained through the formulation with an independent choice of the indices of the operator at each step of the algorithm. Regarding the shuffling heuristic, numerous studies explore methods suitable for addressing classic finite-sum minimization problems. In the work \cite{mishchenko2020random}, the authors examined a classic \textsc{SGD} algorithm, and by introducing a new concept of variance specific to
RR/SO, they matched the lower bounds in such scenarios. In the work \cite{malinovsky2023random}, the \textsc{SVRG} method with RR was considered. The authors actively utilized results from the work \cite{mishchenko2020random} and obtained improved rates. Additionally, there is a set of studies that considered methods incorporating the VR technique in the shuffling setup \cite{huang2021improved,mokhtari2018surpassing,ying2020variance}. However, there are currently no papers that employ the shuffling setting to solve variational inequalities. We aim to fill this gap.

\textbf{Contributions.}
The main results can be summarized as follows.\\
$\bullet$ \textit{Considered algorithms.} We consider two algorithms: \textsc{Extragradient} \cite{juditsky2011solving} and \textsc{Extragradient} that incorporate variance reduction \cite{alacaoglu2022stochastic}, utilizing shuffling heuristics instead of independent index selection.\\
$\bullet$ \textit{Novel approach to proof.} Since shuffling methods lack the property of unbiasedness of stochastic operators, it is essential to propose new approaches to demonstrate convergence. In this paper, we present a technique that enables us to "return" to the starting point of an epoch in which the property of unbiasedness is maintained.\\
$\bullet$ \textit{Convergence estimates.} We provide the first theoretical convergence rates for shuffling methods applied to the finite-sum variational inequality problem. Our comprehensive analysis establishes upper bounds on convergence rates, extending the theoretical framework to encompass this important class of algorithms. In the case of \textsc{Extragradient}, our estimate on the linear term coincides with that for the method without shuffling. In the case of \textsc{Extragradient} with VR, we are the first to obtain a linear convergence estimate for methods with shuffling in the VI problem.\\
$\bullet$ \textit{Experiments.} We conduct comprehensive experiments that emphasize the superiority of shuffling over the random index selection heuristic. We consider two classic practical applications: image denoising and adversarial training.

\vspace{-2mm}
\section{Setup}
\vspace{-2mm}

\textbf{Assumptions.} We present a list of assumptions under which we derive the main statements.
\begin{assumption}\label{as:lipschitz}
    Each operator $F_i$ is $L$-Lipschitz, i.e., it satisfies $\|F_i(z_1) - F_i(z_2)\| \leq L\|z_1 - z_2\|$ for any $z_1, z_2 \in Z$.
\end{assumption}
\begin{assumption}\label{as:monotone}
    Each operator $F_i$ is $\mu$-strongly monotone, i.e., it satisfies $\langle F_i(z_1) - F_i(z_2), z_1 - z_2\rangle \geqslant \mu\|z_1 - z_2\|^2$ ~for any $z_1, z_2 \in Z$.  
\end{assumption}
\begin{assumption}\label{as:bound}
    Each stochastic operator $F_i$ and full operator $F$ is bounded at the point of the solution $z^*\in\text{dom}~ g$, i.e. $\mathbb E \|F_i(z^*)\|^2 \leqslant \sigma_*^2, \|F(z^*)\|^2 \leqslant \sigma_*^2$.  
\end{assumption}

\textbf{Proximal Algorithm.}
Earlier, we provided examples of the application of variational inequalities (Examples \ref{ex:convopt}, \ref{ex:convconcsaddle}). In many optimization problems, particularly in machine learning and signal processing, we often encounter the need to minimize a function of the same form, i.e., decompose it into two parts: a smooth differentiable function \( f : \mathbb{R}^n \rightarrow \mathbb{R} \) and a possibly non-smooth function \( g : \mathbb{R}^n \rightarrow \mathbb{R} \). To solve this problem, we utilize the proximal gradient method. The core idea is to iteratively update the solution by combining \textsc{Gradient descent} on the smooth part \( f \) and the proximal operator for the possibly non-smooth part \( g \). We also assume that $g$ is proximal friendly, meaning the solution of the minimization problem on $g$ is achieved at minimal cost. The proximal operator of the function \( g \) at a point \( x \) is defined as:
\begin{equation*}
\text{prox}_{g}(z) = \arg\min_{y \in \mathbb{R}^n} \left\{ g(y) + \frac{1}{2} \| y - z \|^2 \right\},
\end{equation*}
where \( \| \cdot \| \) denotes the Euclidean norm.
Using the proximal operator, the update step for solving the optimization problem can be expressed as
\begin{equation*}
z^{t+1} = \text{prox}_{\alpha_t g}\left( z^t - \alpha_t \nabla f(z^t) \right).
\end{equation*}
For us, the proximal operator plays a role, since \eqref{eq:vi_setting} also uses a regularizer.

\section{Algorithms and Convergence Analysis}
\subsection{Extragradient}

The setting of shuffling lies in the fact that we do not choose the stochastic operator independently at each step of the method. Instead, we permute the sequence of indices and, at each iteration of the algorithm, we select the operator according to the new sequence. In this work, we focus on the Random Reshuffling and Shuffle Once techniques and provide appropriate \textsc{Extragradient} methods (Algorithms \ref{alg:rrextragrad}, \ref{alg:soextragrad}).  
%\vspace{-6mm}

\begin{minipage}[t]{0.48\linewidth}
    \begin{algorithm}{\textsc{RR Extragradient}}
    \label{alg:rrextragrad}
    \begin{algorithmic}[1]
    \State \textbf{Input:} Starting point $z^0_0\in\mathbb{R}^d$
    \State \textbf{Parameter:} Stepsize $\gamma$
    \For{$s = 0, 1, 2, \ldots, S-1$}
        \State \textcolor{orange}{Generate a permutation $\pi_0, \pi_1, \ldots, \pi_{n-1}$ of sequence $\{1, 2, \ldots, n\}$}
        \For{$t = 0, 1, 2, \ldots, n-1$}  
            \State $z^{t + \frac{1}{2}}_s = \text{prox}_{\gamma g}\left(z^t_s - \gamma F_{\pi_s^t}(z^t_s)\right)$
            \State $z^{t + 1}_s = \text{prox}_{\gamma g}\left(z^t_s - \gamma F_{\pi_s^t}(z^{t + \frac{1}{2}}_s)\right)$
        \EndFor
        \State $z_s^n = z_{s+1}^0$
    \EndFor
    \State \textbf{Output:} $z^n_S$
    \end{algorithmic}
    \end{algorithm}
\end{minipage}
\hfill
\begin{minipage}[t]{0.48\linewidth}
    \begin{algorithm}{\textsc{SO Extragradient}}
    \label{alg:soextragrad}
    \begin{algorithmic}[1]
    \State \textbf{Input:} Starting point $z^0_0\in\mathbb{R}^d$
    \State \textbf{Parameter:} Stepsize $\gamma$
    \State \textcolor{orange}{Generate a permutation $\pi_0, \pi_1, \ldots, \pi_{n-1}$ of sequence $\{1, 2, \ldots, n\}$}
    \For{$s = 0, 1, 2, \ldots, S-1$}
        \For{$t = 0, 1, 2, \ldots, n-1$}  
            \State $z^{t + \frac{1}{2}}_s = \text{prox}_{\gamma g}\left(z^t_s - \gamma F_{\pi_s^t}(z^t_s)\right)$
            \State $z^{t + 1}_s = \text{prox}_{\gamma g}\left(z^t_s - \gamma F_{\pi_s^t}(z^{t + \frac{1}{2}}_s)\right)$
        \EndFor
        \State $z_s^n = z_{s+1}^0$
    \EndFor
    \State \textbf{Output:} $z^n_S$
    \end{algorithmic}
    \end{algorithm}
    \vspace{2mm}
\end{minipage}

The analysis of shuffling methods has specific details. The key difference between shuffling and independent choice is that shuffling methods lack one essential feature: the unbiasedness of stochastic operators:
\begin{equation*}
    \mathbb E_{\pi_s^t} \left[F_{\pi_s^t}(z_s^t)\right] \neq \frac{1}{n}\sum\nolimits_{i = 1}^n F_{\pi_s^i}(z_s^t) = F(z_s^t).
\end{equation*}

\vspace{-2mm}
This restriction leads us to a more complex analysis and the use of non-standard techniques to prove the convergence of the shuffling methods. Nevertheless, at two points -- $z_s^0$ and $z^*$, this equality holds. Indeed, the point $z_s^0$ is the initial point of the epoch, where we choose one random index from $n$. Additionally, the point $z^*$ does not depend on $t$. Thus, we can "go back" to the beginning of the epoch and leverage the unbiased operators. This technique is interesting not only in relation to shuffling methods. For example, it is applicable to methods that utilize Markov chains to select indices, as there is no unbiased property except at the correlation point of the chain. This is the key point of our analysis, and now, having established this, we present the main result of this section.

\begin{theorem}\label{th:eg}
    Suppose Assumptions \ref{as:lipschitz}, \ref{as:monotone}, \ref{as:bound} hold. Then for Algorithms \ref{alg:rrextragrad}, \ref{alg:soextragrad} with $\gamma\leqslant\min\left\{\frac{1}{2\mu n}, \frac{1}{6L}\right\}$ after $S$ epochs,
    \vspace{-2mm}
    \begin{equation*}
        \|z_S^n - z^*\|^2 \leqslant (1 - \frac{\gamma\mu}{2})^{Sn} \|z^0_0 - z^*\|^2 +  \frac{256\gamma n^2\sigma^2_*}{\mu}.
    \end{equation*}
\end{theorem}

\begin{corollary}\label{cor:eg}
    Suppose Assumptions \ref{as:lipschitz}, \ref{as:monotone}, \ref{as:bound} hold. Then Algorithms \ref{alg:rrextragrad}, \ref{alg:soextragrad} with $\gamma \leqslant \min\left\{\frac{1}{2\mu n}, \frac{1}{6L}, \frac{2\log\left(\max\left\{2, \frac{\mu^2\|z_0^0 - z^*\|^2 T}{512 n^2\sigma^2_*}\right\}\right)}{\mu T}\right\}$, to reach $\varepsilon$-accuracy, where $\varepsilon \sim \|z_S^n - z^*\|^2$, needs
    \begin{equation*}
    \mathcal{\widetilde{O}}\left(\left(n + \frac{L}{\mu}\right)\log\left(\frac{1}{\varepsilon}\right) + \frac{n^2\sigma^2_*}{\mu^2\varepsilon}\right) ~~\text{iterations and oracle calls.}
    \end{equation*}
\end{corollary}

\begin{remark}\label{rem:eg}
    We can transform the obtained estimation for the case of monotone stochastic operators. To achieve this, we utilize a regularization trick with $\mu\sim \frac{\varepsilon}{D}$. In particular, solving the problem with the operator $\hat{F}(z) = F(z) + \mu(z-z_0^0)$ at accuracy $\frac{\varepsilon}{2}$ allows us to solve the problem \eqref{eq:vi_setting} at accuracy $\varepsilon$, resulting in $\mathcal{\widetilde{O}}\left(n + \frac{L}{\varepsilon} + \frac{n^2}{\varepsilon^3}\right)$ for iteration and oracle complexity. This represents convergence in argument, which differs from the classic form.
\end{remark}

Let us explain the result of the theorem. The form of the estimate is classic and appears in all stochastic methods for strongly convex minimization \cite{moulines2011non,stich2019unified} and strongly monotone VIs \cite{beznosikov2020distributed,mishchenko2020revisiting}. We compare it with the results of related works. Our method is based on \textsc{REG} \cite{mishchenko2020revisiting}. In this work, authors obtain $\mathcal{\widetilde{O}}\left(\frac{L}{\mu} + \frac{1}{\mu^2\varepsilon}\right)$ oracle complexity. Thus, our result represents a significant advancement in shuffling theory. Although there is a degradation in $n$ in the sublinear term, the estimation on the linear term coincides with that in the classic setting using an independent choice of stochastic operators. We also compare the result with the work \cite{juditsky2011solving}. The authors obtain $\mathcal{O}\left(\frac{L}{\sqrt{\varepsilon}} + \frac{1}{\varepsilon}\right)$. However, uniform bounds on the variance were required in this work, while we bound the variance only at the optimum. Note that, according to current theory, shuffling methods are no more effective than methods with independent sampling for classic minimization problems. Indeed, classic \textsc{SGD} rate in the non-convex case is $\mathcal{O}\left(\frac{L}{\sqrt{\varepsilon}} + \frac{L}{\varepsilon}\right)$. Therefore, prior works on shuffling in the minimization setup deliver $\mathcal{O}\left(\frac{nL}{\sqrt{\varepsilon}} + \frac{nL}{\varepsilon^{\nicefrac{3}{4}}}\right)$ \cite{mishchenko2020random,mohtashami2022characterizing,lu2022general} and $\mathcal{O}\left(\frac{L}{\sqrt{\varepsilon}} + \frac{nL}{\varepsilon^{\nicefrac{3}{4}}}\right)$ \cite{koloskova2023convergence} rates. Thus, despite the improvement in the asymptotics of the variance term, it deteriorates the theoretical estimate by a factor of $n$.

Let us focus on the second term in the estimation.  In general, the sublinear term with $\sigma_*^2$ is not improved. However, for the finite-sum problem, this term can be eliminated by employing additional techniques, such as variance reduction.

\subsection{Extragradient with Variance Reduction}

%\begin{wrapfigure}{r}{\textwidth}
%\begin{minipage}{\textwidth}
\begin{algorithm}{RR/SO \textsc{Extragradient} with variance reduction}
\label{alg:proxextragradvr}
\begin{algorithmic}[1]
\State \textbf{Input:} \textbf{Parameters:} $z_0^0, \omega_0^0$
\State \textbf{Parameter:} Stepsize $\gamma$, $\alpha \in (0, 1)$ 
\State \textcolor{orange}{Generate a permutation $\pi_0, \pi_1, \ldots, \pi_{n-1}$ of sequence $\{1, 2, \ldots, n\}$} \quad\quad//\quad SO heuristic
\For{$s = 0, 1, \ldots$}
    \State \textcolor{orange}{Generate a permutation $\pi_0, \pi_1, \ldots, \pi_{n-1}$ of sequence $\{1, 2, \ldots, n\}$} ~~//\quad RR heuristic
    \For{$t = 0, 1, \ldots, n-1$}
        \State $\overline{z}_s^t = \alpha z_s^t + (1 - \alpha) \omega_s^t$
        \State$z_s^{t+\nicefrac{1}{2}} = \text{prox}_{\gamma g} \left(\overline{z}_s^t - \gamma F\left(\omega_s^t\right)\right)$
        \State $\hat{F}(z_s^{t + \nicefrac{1}{2}}) = F_{\pi_s^t} (z_s^{t + \nicefrac{1}{2}}) - F_{\pi_s^t} (\omega_s^t) + F(\omega_s^t)$
        \State$z_s^{t+1} = \text{prox}_{\gamma g} \left(\overline{z}_s^t - \gamma \hat{F} \left(z_s^{t+\nicefrac{1}{2}}\right)\right)$
        \State\label{alg3:line11} $\omega_s^{t+1} = \begin{cases}
            z_s^t, \quad &\text{with probability}\quad p\\
            \omega_s^t \quad &\text{with probability}\quad 1-p
        \end{cases}$
    \EndFor
    \State $z_{s+1}^0 = z_s^n$
    \State $\omega_{s+1}^0 = \omega_s^n$
\EndFor
\State \textbf{Output:} $z^n_S$
\end{algorithmic}
\end{algorithm}
%\end{minipage}
%\end{wrapfigure}

Now we utilize the variance reduction technique, which enhances the convergence of the algorithms by diminishing the impact of random fluctuations. This approach was not employed in the previously presented algorithms. We introduce a variant of the RR/SO \textsc{Extragradient} with a variance reduction algorithm (Algorithm \ref{alg:proxextragradvr}) and provide the convergence results for this method. In the work \cite{malinovsky2023random}, where shuffling is investigated in variance reduction methods, the authors utilize a more classic version and compute $F(\omega_s^t)$ at the beginning of each epoch. We consider an alternative and compute this full operator randomly with probability $p$. We set $p = \frac{1}{n}$, which implies that, on average, the full operator is updated once per epoch. Note that this choice does not increase the oracle complexity.

\begin{theorem}\label{th:proxegvr}
    Suppose that Assumptions \ref{as:lipschitz}, \ref{as:monotone} hold. Then for Algorithm \ref{alg:proxextragradvr} with $\gamma \leqslant\frac{(1-\alpha)\mu}{6L^2}$ and $p = \frac{1}{n}$ after $T$ iterations,
    \vspace{-1mm}
    \begin{equation*}
        V_S^n \leqslant \left(1 - \frac{\gamma\mu}{4}\right)^T V_0^0,
    \end{equation*}
    where $V_s^t =\mathbb E \|z_s^t - z^*\|^2 + \mathbb E\|\omega_s^t - z^*\|^2$.
\end{theorem}

\begin{corollary}\label{cor:proxegvr}
    Suppose that Assumptions \ref{as:lipschitz}, \ref{as:monotone} hold. Then Algorithm \ref{alg:proxextragradvr} with $\gamma \leqslant\frac{(1-\alpha)\mu}{6L^2}$ and $p = \frac{1}{n}$, to reach $\varepsilon$-accuracy, where $\varepsilon \sim V_S^n$, needs
    \vspace{-1mm}
    \begin{align*}
    &\mathcal{O}\left(n\frac{L^2}{\mu^2}\log\left(\frac{1}{\varepsilon}\right)\right) ~~\text{iterations and oracle calls,}
    \end{align*}
    where $V_s^t =\mathbb E \|z_s^t - z^*\|^2 + \mathbb E\|\omega_s^t - z^*\|^2$.
\end{corollary}

\begin{remark}\label{rem:egvr}
    Similarly to Remark \ref{rem:eg}, we can apply our result in the monotone case through the regularization trick and obtain $\mathcal{\widetilde{O}}\bigl(n\frac{L^2}{\varepsilon^2}\bigr)$.
\end{remark}

We remove the variance that arose in Theorem \ref{th:eg} and obtain linear convergence. Even though we obtain worse estimates than those in works that also use the variance reduction technique, such as \cite{alacaoglu2022stochastic,alacaoglu2021forward,chavdarova2019reducing,palaniappan2016stochastic} (see Table \ref{tab:methodscomp}), there is a distinct explanation for this. According to current theory, methods with the shuffling heuristic are inferior to methods with independent sampling for the variance reduction methods \cite{malinovsky2023random}. In the prior  works, where the minimization of the strongly-convex objective with shuffling was considered, the following rates are obtained: $\mathcal{\widetilde{O}}\left(n^2\frac{L^2}{\mu^2}\right)$ \cite{gurbuzbalaban2017convergence}, $\mathcal{\widetilde{O}}\left(n\frac{L^2}{\mu^2}\right)$ \cite{ying2020variance}, $\mathcal{\widetilde{O}}\left(n\frac{L^{\nicefrac{3}{2}}}{\mu^{\nicefrac{3}{2}}}\right)$ \cite{malinovsky2023random}. At the same time, in classic minimization, variance reduction methods provide $\mathcal{\widetilde{O}}\left(n\frac{L}{\mu}\right)$ rates \cite{gorbunov2020unified}. Thus, it remains an open question whether it is possible to obtain theoretical convergence estimates for methods using shuffling heuristics that are equivalent to those for methods with independent index selection for the VI problem. Additionally, in the course of the work, no theoretical differences are revealed in the SO and RR techniques concerning the problem \eqref{eq:vi_setting}.

\section{Experiments}
\label{sec:experiments}
In this section, we evaluate the proposed algorithms to demonstrate their practical applications by conducting experiments in two cases: image denoising and adversarial training.

\subsection{Image Denoising} 

To formulate the image denoising problem \cite{chambolle2011first}, we consider the classic saddle point problem as presented in Example \ref{ex:convconcsaddle}:
$$
\min_{x \in \mathcal{X}} \max_{y \in \mathcal{Y}} \left[\langle Kx, y \rangle + G_1(x) - G_2(y)\right],
$$
where regularizers \(G_1\) and \(G_2\) are proper convex lower semicontinuous functions, and \(K\) is a continuous linear operator. To proceed with image denoising, we consider $g$ as a given noisy image and $u$ as the solution we seek. We use the Cartesian grid with the step $h: \{(i\cdot h, j\cdot h)\}$. Thus, specifically for image denoising, we consider:
\vspace{-1mm}
$$
\min_{u \in \mathcal{X}} \max_{p \in \mathcal{Y}} \left[\langle \nabla u, p \rangle_{\mathcal{Y}} + \nicefrac{\lambda}{2} \|u - g\|^2_2 - \delta_{P}(p)\right],
\vspace{-2mm}
$$
where $p$ is a dual variable, \(\delta_{P}(p)\) is the indicator function of the set \(P\) defined as: $ P = \{ p \in \mathcal{Y} : \| p(x) \| \leq 1 \}.$ The indicator function \(\delta_P(p)\) is defined as zero if \(p\) belongs to the set \(P\), and infinity otherwise. We define operator $\nabla u$ as the difference between neighboring pixels in the grid, both horizontally and vertically, normalized by the step of the grid $h$. This formulation represents a saddle point problem, where we seek to minimize the first term with respect to \(u\) while simultaneously maximizing the second term with respect to \(p\). Using duality, we can write the final formulation of the problem as
\begin{equation} \label{problem_denoising}
\min_{u \in \mathcal{X}} \max_{p \in \mathcal{Y}} \left[-\langle u, \text{div} ~p \rangle_{\mathcal{X}} + \nicefrac{\lambda}{2} \|u - g\|^2_2 - \delta_{P}(p)\right].
\end{equation}
\vspace{-4mm}

To bring the problem to the form of a finite sum \eqref{eq:finite-sum}, we divide images into batches of equal squares. We consider two options: batches of size 4 and 8, according to the grid. Since the images are black and white, they are single-channel, which means that each batch is a square matrix with non-negative integers. It is also important to note that when calculating the gradient, the edges of the batch are processed according to the rule of adding a number equal to the nearest neighbor.

We compare RR/SO \textsc{Extragradient} with variance reduction (Algorithm \ref{alg:proxextragradvr}) and \textsc{Extragradient} with variance reduction \cite{alacaoglu2022stochastic}. Analogously, we compare RR/SO \textsc{Extragradient} (Algorithms \ref{alg:rrextragrad}, \ref{alg:soextragrad}) and \textsc{Extragradient} \cite{juditsky2011solving}. We select two images with different levels of additive zero-mean Gaussian noise: \(\sigma = 0.05\) and \(\sigma = 0.1\). Figures \ref{fig:eggirl} and \ref{fig:mpvrgirl} provide a comparison of the proposed methods.
Additional results for all considered methods on another image are presented in Figures \ref{fig:egluvr}, \ref{fig:mpvrluvr} in \hyperref[A]{Appendix A}. 

Comparing the images, it is evident that algorithms incorporating shuffling perform better than those that do not, although the difference in the line graphs is subtle. While the results from the Independent Choice strategy appear sharper, they are also noisier compared to those from algorithms employing shuffling. Thus, utilizing shuffling techniques reduces noise more effectively.

\begin{figure}[h]
    \centering
    \includegraphics[width=0.83\linewidth]{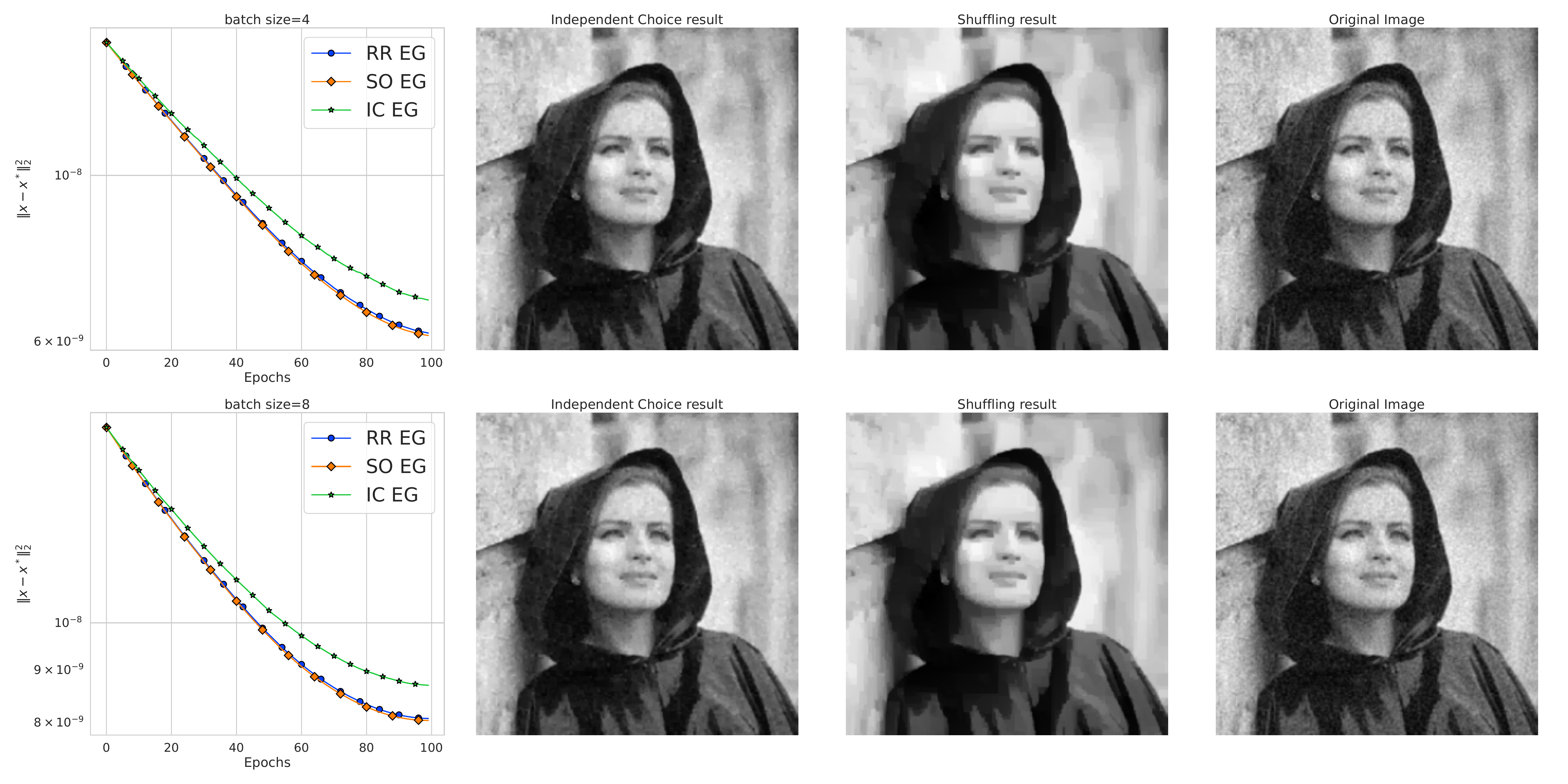}
    \caption{\textsc{Extragradient} convergence on image with $\sigma = 0.05$ on the problem \eqref{problem_denoising}.}
    \label{fig:eggirl}
\end{figure}

\begin{figure}[h]
    \centering
    \includegraphics[width=0.83\linewidth]{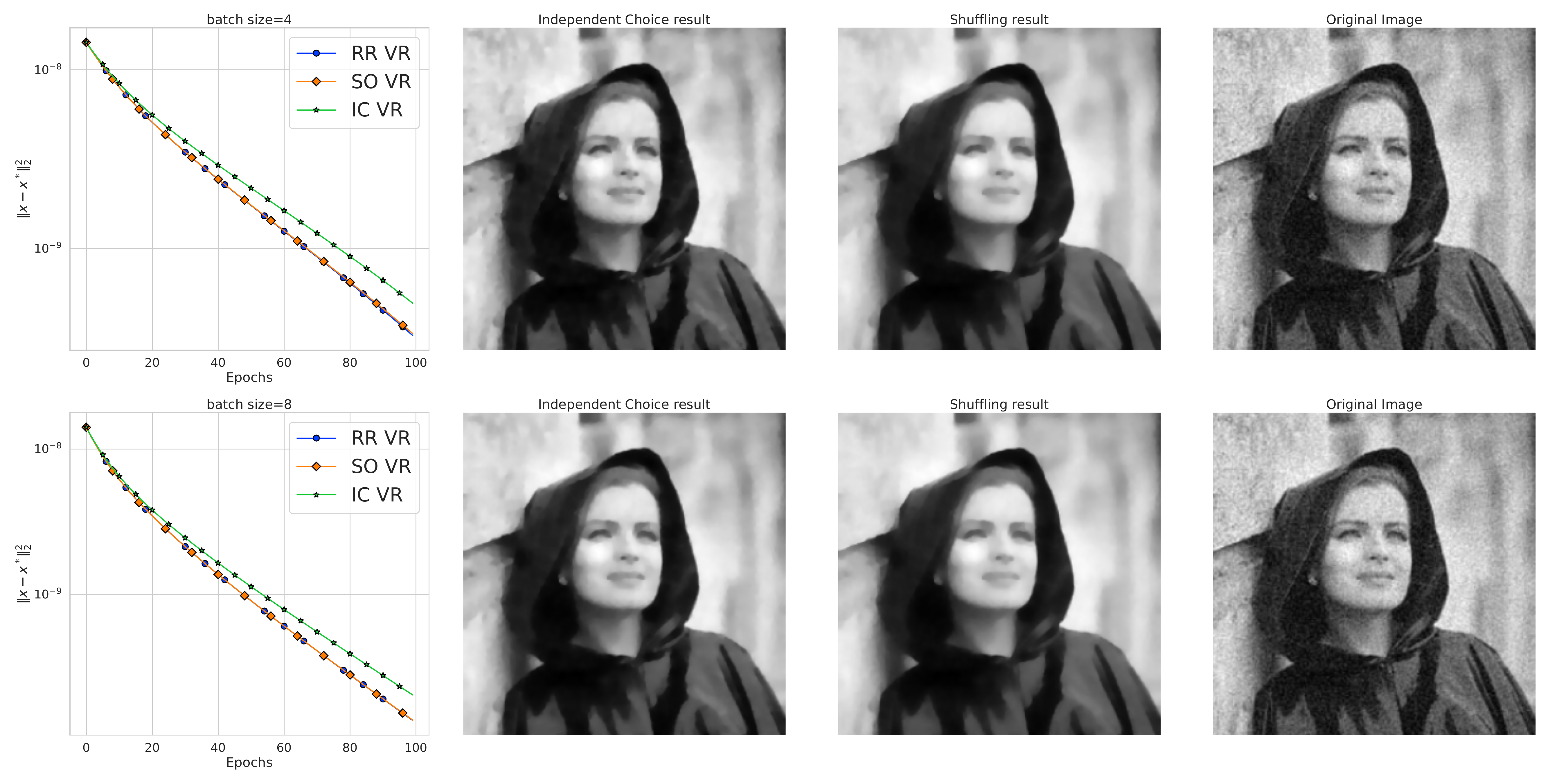}
    \caption{\textsc{Extragradient} with VR convergence on image with $\sigma = 0.05$ on the problem \eqref{problem_denoising}.}
    \label{fig:mpvrgirl}
\end{figure}

% \vspace{-2mm}
% \begin{center}
% \includegraphics[width=0.83\linewidth]{vr_ref_eg_second.pdf}
% \vspace{-4mm}
% \captionof{figure}{\textsc{Extragradient} convergence on image with $\sigma = 0.05$ on the problem \eqref{problem_denoising}.}
% \label{fig:eggirl}
% \end{center}
% \vspace{-4mm}
% \begin{center}
% \includegraphics[width=0.83\linewidth]{vr_ref_second.pdf}
% \vspace{-4mm}
% \captionof{figure}{\textsc{Extragradient} with VR convergence on image with $\sigma = 0.05$ on the problem \eqref{problem_denoising}.}
% \label{fig:mpvrgirl}
% \end{center}
%\vspace{-1mm}
%\vspace{-2mm}

\paragraph{Hyperparameter details.}
We configure the training with a fixed learning rate $\gamma = 5 \times 10^{-3}$ and a batch size of $8$. The probability parameter $p = \frac{1}{n}$ is determined based on the effective dataset size. However, for computational efficiency, we adapt it to the mini-batch context. Training is conducted with random state $50$.

\subsection{Adversarial Training}

Next, we address an adversarial training problem. We formulate it in the following manner: 

\begin{equation} \label{problem_adversarial}
\min_{w \in \mathbb{R}^d} \max_{\|r_i\| \leqslant D} \left[\frac{1}{2N} \sum \limits_{i=1}^N \left( w^T\left(x_i+r_i\right) - y_i\right) ^2 + \frac{\lambda}{2} \|w\|^2 -\frac{\beta}{2}\|r\|^2\right],
\end{equation}
where the samples correspond to features $x_i$ and targets $y_i$. We evaluate this issue across several datasets: \texttt{mushrooms}, \texttt{a9a}, and \texttt{w8a}, sourced from the \textsc{LIBSVM} library \cite{chang2011libsvm}. A brief description of these datasets is provided in Table \ref{tab:dataset-summary}, \hyperref[A]{Appendix A}. The results are presented in Figure \ref{fig:adversarial}.

\begin{figure}[htbp]
    \centering
    \begin{subfigure}[b]{0.3\textwidth}
        \centering
        \includegraphics[width=\textwidth]{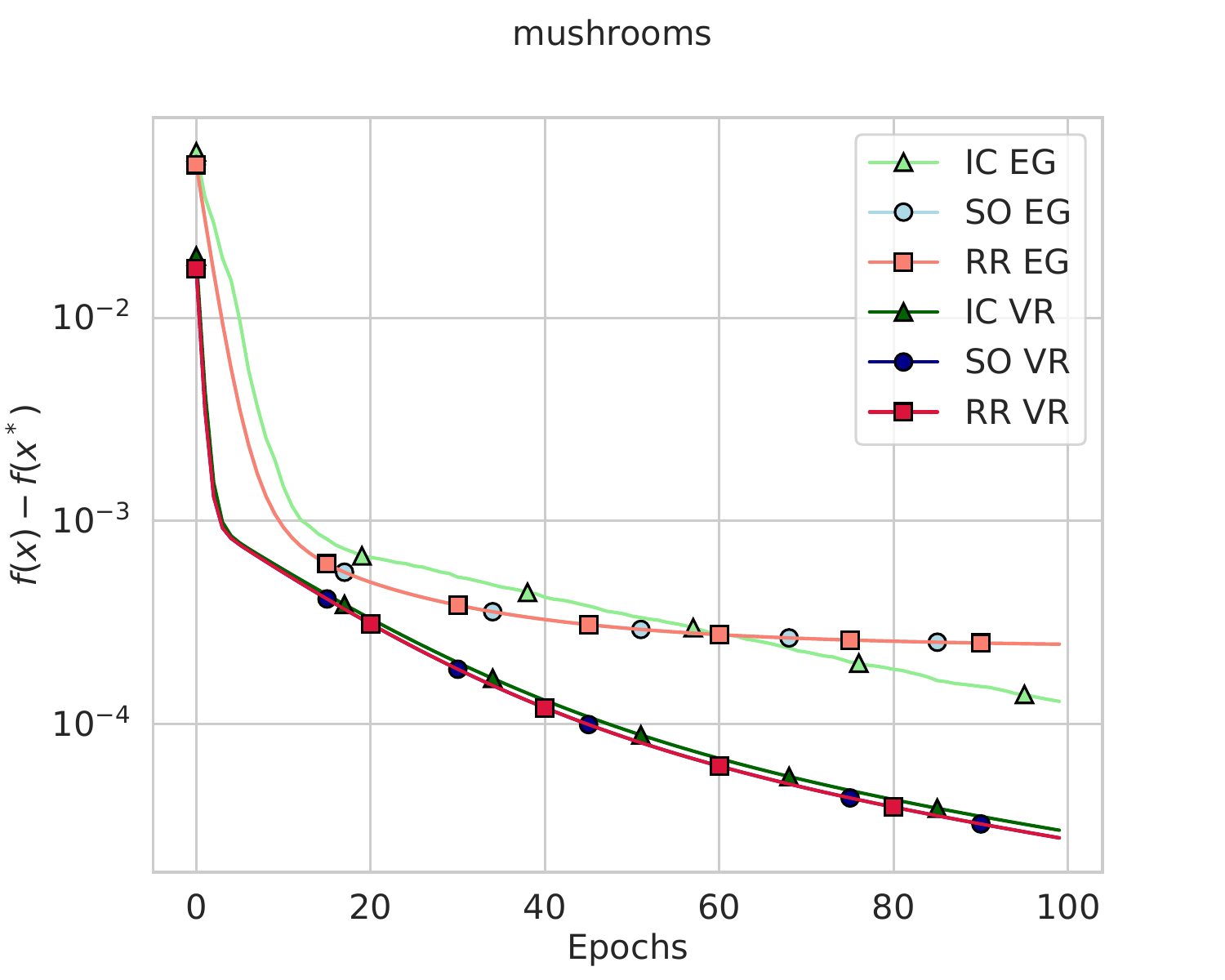}
        \caption{\texttt{mushrooms}}
    \end{subfigure}
    \hfill
    \begin{subfigure}[b]{0.3\textwidth}
        \centering
        \includegraphics[width=\textwidth]{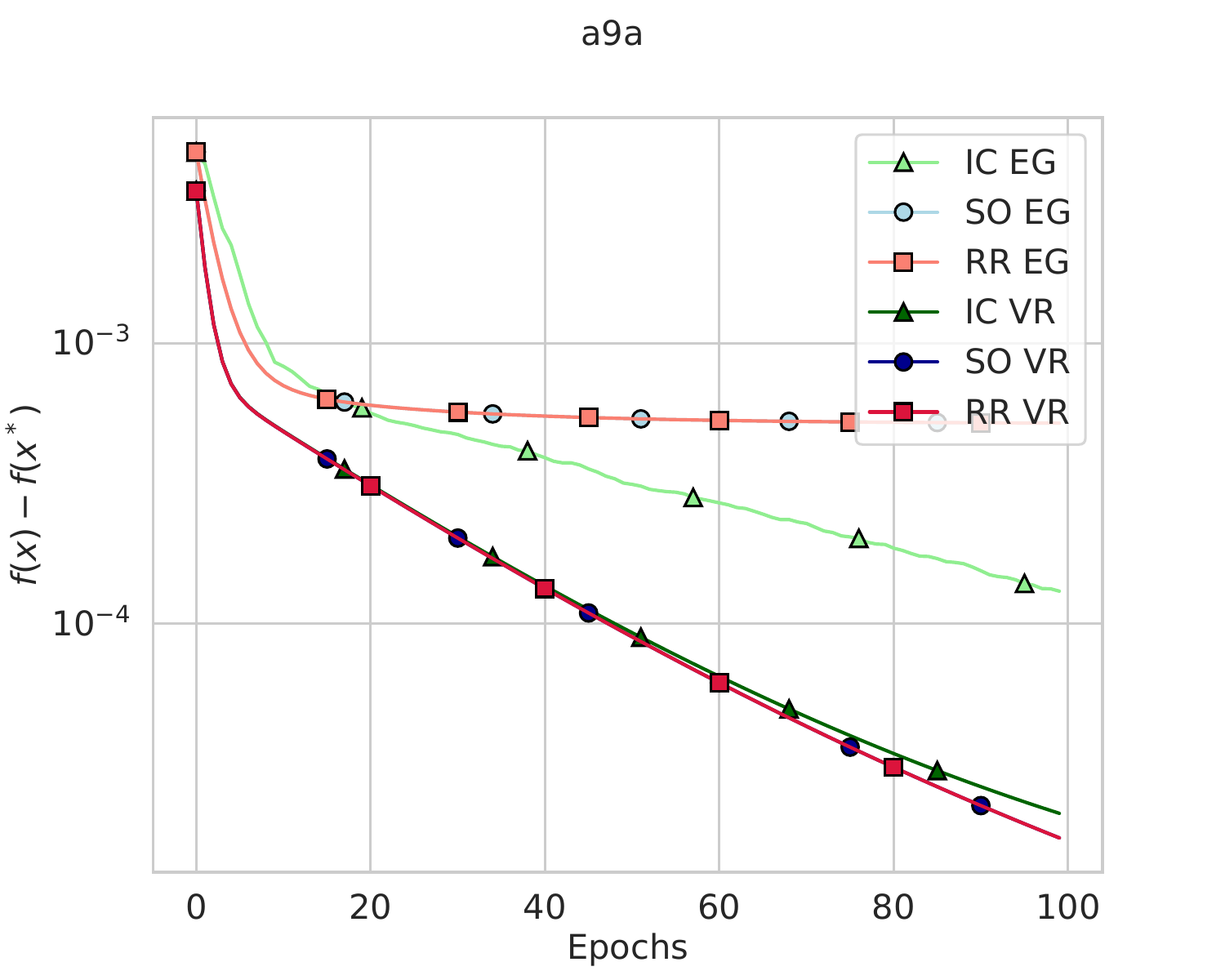}
        \caption{\texttt{a9a}}
    \end{subfigure}
    \hfill
    \begin{subfigure}[b]{0.3\textwidth}
        \centering
        \includegraphics[width=\textwidth]{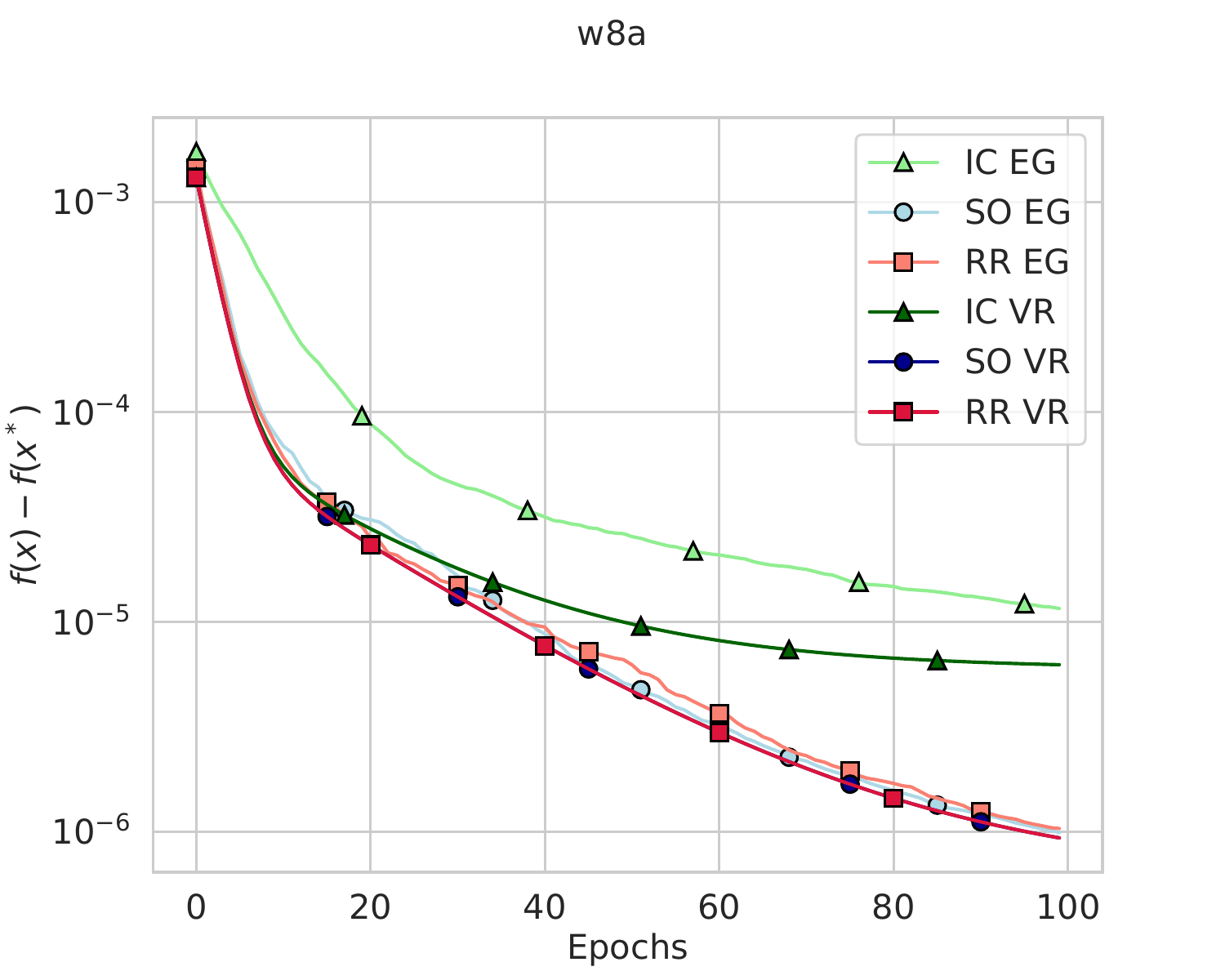}
        \caption{\texttt{w8a}}
    \end{subfigure}
    \caption{Extragradient with and without VR compared using various shuffling heuristics on the datasets shown above for the problem \eqref{problem_adversarial}.}
    \label{fig:adversarial}
\end{figure}

As shown in the plots, \textsc{Extragradient} and \textsc{Extragradient} with the VR algorithms using the independent choice of indices demonstrate worse performance compared to those using shuffling methods. In this series of experiments, RR exhibits better performance than other shuffling methods and significantly outperforms the non-shuffled versions.

\textbf{Hyperparameter Details.}
We configure the training with a fixed learning rate $\gamma = 0.01$ and a batch size of $4,4,16$ on  \texttt{mushrooms}, \texttt{a9a} and \texttt{w8a} datasets, respectively. The probability parameter $p = \frac{1}{n}$ is determined based on the effective dataset size; however, for computational efficiency, we adapt it to the mini-batch context. Training is conducted with random state $50$.

\section{Conclusion}

In this paper, we consider stochastic algorithms for solving the problem of variational inequalities. Specifically, we explore the influence of shuffling heuristics on stochastic methods. For the first time, a theoretical analysis of two versions of the \textsc{Extragradient} algorithm, both with and without variance reduction involving shuffling, is presented. Empirical results on image denoising and adversarial training tasks confirm the applicability of the methods in practice. Nevertheless, the convergence rate of the \textsc{Extragradient} algorithm with variance reduction is not optimal (see Corollary \ref{cor:proxegvr}). Therefore, its improvement can be considered as a future work.

\end{mainpart}

\begin{appendixpart}

\section{Additional Experiments} \label{A}

In this section, we present additional experiments that have been performed (Figures \ref{fig:egluvr}, \ref{fig:mpvrluvr}). Similar to the previous experiments, we observe a consistent pattern: methods incorporating shuffling techniques outperform those without shuffling. These results further confirm the effectiveness of shuffling techniques in addressing the denoising problem on another image with higher $\sigma$.
\begin{center}
\includegraphics[width=0.9\linewidth]{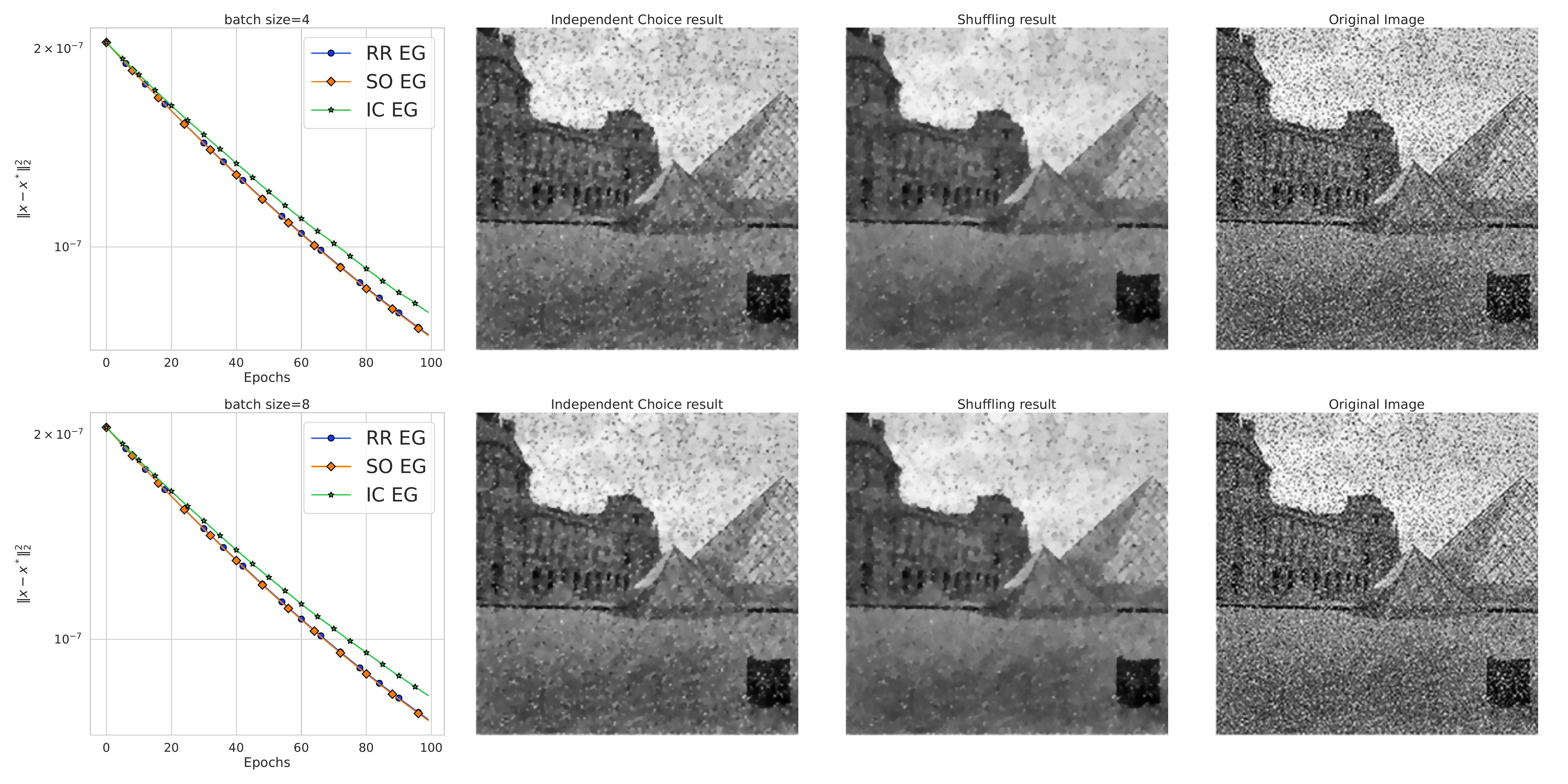}
\captionof{figure}{\textsc{Extragradient} convergence on image with  $\sigma = 0.1$ on the problem \eqref{problem_denoising}.}
\label{fig:egluvr}
\end{center}

\begin{table}[!htbp]
    \centering
    \begin{tabularx}{\textwidth}{lXXX}
        \toprule
        \textbf{Name} & \begin{tabular}{@{}l@{}}
\textbf{Number of} \\
\textbf{Instances}
\end{tabular} & \begin{tabular}{@{}l@{}}
\textbf{Number of} \\
\textbf{Features}
\end{tabular} & \begin{tabular}{@{}l@{}}
\textbf{Number of} \\
\textbf{Classes}
\end{tabular} \\
        \midrule
        \texttt{mushrooms} & 8,124 & 112 & 2 \\
        \texttt{a9a} & 32,561 & 123 & 2 \\
        \texttt{w8a} & 49,749 & 300 & 2 \\
        \bottomrule
    \end{tabularx}
    \caption{Summary of Datasets}
    \label{tab:dataset-summary}
\end{table}

\begin{center}
\includegraphics[width=0.9\linewidth]{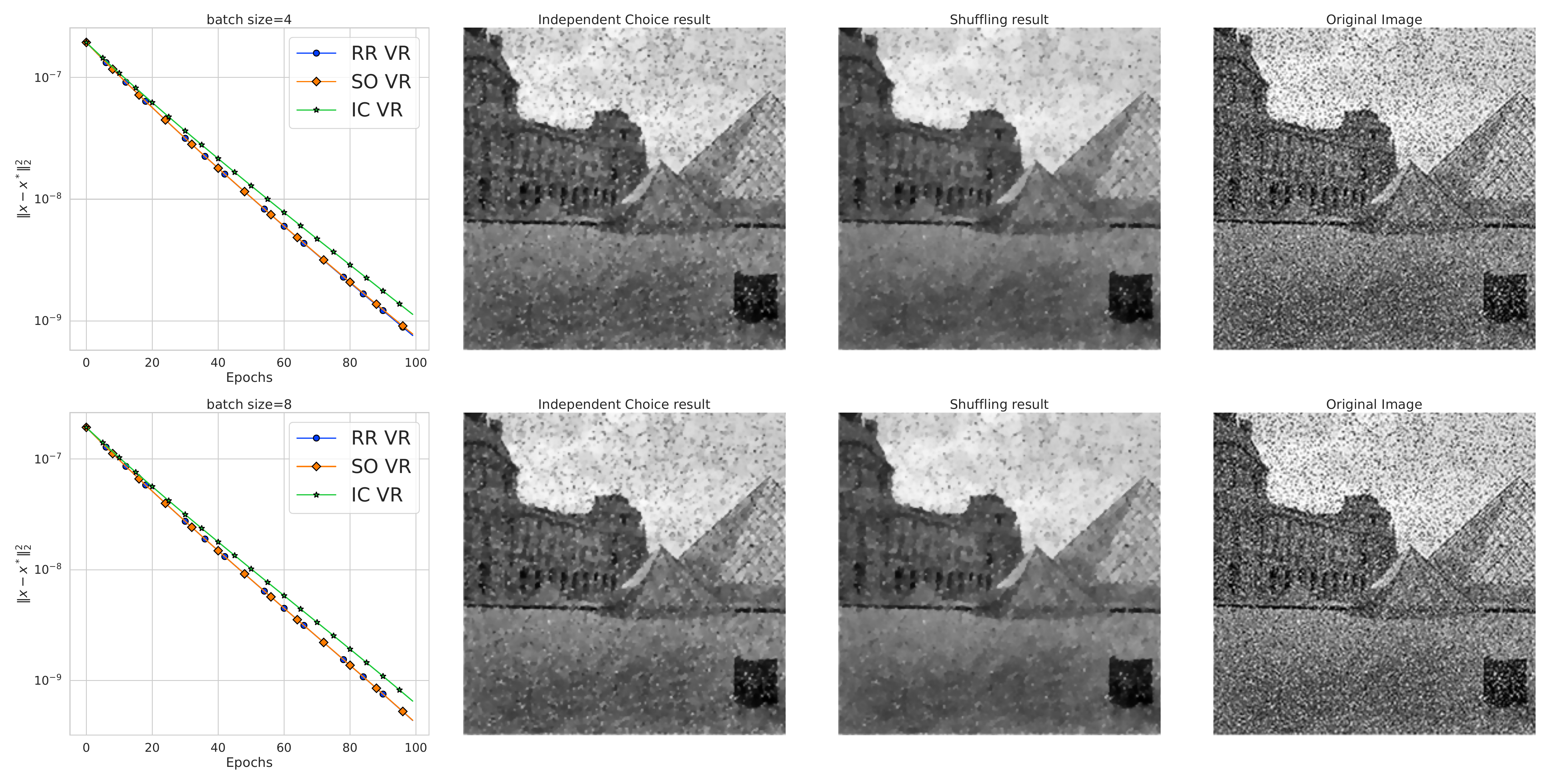}
\captionof{figure}{\textsc{Extragradient} with VR convergence on image with $\sigma = 0.1$ on the problem \eqref{problem_denoising}.}
\label{fig:mpvrluvr}
\end{center}

The datasets used for the experiments on the adversarial training include \texttt{mushrooms}, \texttt{a9a}, and \texttt{w8a}. These datasets vary in size and complexity (see Table \ref{tab:dataset-summary} for details), providing a comprehensive evaluation of our proposed algorithms in the context of adversarial training.

\section{Basic Inequalities}\label{B}

For all vectors \( x, y, \{x_i\}_{i=1}^n \) in \( \mathbb{R}^d \) with a positive scalar \( \alpha \), the following holds:
\begin{align}
\label{bi:1} \tag{Scalar} \langle x, y \rangle ~ & \leqslant ~ \frac{\|x\|^2}{2\alpha} + \frac{\alpha \|y\|^2}{2}, \\
\label{bi:quadr} \tag{Norm} 2\langle x, y \rangle ~ & = ~ \|x+y\|^2 - \|x\|^2 - \|y\|^2, \\
\label{bi:quadrineq} \tag{CS} -2\|x\|^2 ~ & \leqslant ~ -\|x+y\|^2 + 2\|y\|^2, \\
\label{bi:CauchySchwarz} \tag{Sum}  \left\|\sum_{i=1}^{n} x_i\right\|^2 & \leqslant  ~n \sum_{i=1}^{n} \|x_i\|^2,\\
\label{bi:Prox} \tag{Prox}  \left\|\text{prox}_{\gamma g} (x) - \text{prox}_{\gamma g} (y)\right\|^2 & \leqslant  \left\|x-y\right\|^2.
\end{align}

\section{Extragradient}\label{C}
\textbf{Theorem \ref{th:eg}.}
\textit{Suppose Assumptions \ref{as:lipschitz}, \ref{as:monotone} hold. Then for Algorithms \ref{alg:rrextragrad}, \ref{alg:soextragrad} with $\gamma\leqslant \min \left\{\frac{1}{2\mu n}, \frac{1}{6L} \right\}$ after $S$ epochs,}
    \begin{equation*}
        \|z_S^n - z^*\|^2 \leqslant (1 - \frac{\gamma\mu}{2})^{Sn} \|z^0_0 - z^*\|^2 +  \frac{256\gamma n^2\sigma^2_*}{\mu}.
    \end{equation*}

    \begin{proof}

         We begin with the standard prox-inequality:
        \begin{equation}\label{ineq:prox}
            \hat{z} = \text{prox}_{g}(z) \Longleftrightarrow \langle \hat{z} - z, u - \hat{z}\rangle \geqslant g(\hat{z}) - g(u), \quad\forall u\in \mathcal{Z}.
        \end{equation}
        Substituting both steps of Algorithm \ref{alg:rrextragrad} (Algorithm \ref{alg:soextragrad}) into \eqref{ineq:prox}, we derive:
        \begin{align*}
            \langle z_s^{t+1} - z_s^t + \gamma F_{\pi_s^t}(z_s^{t+\nicefrac{1}{2}}), z^* - z_s^{t+1}\rangle &\geqslant \gamma(g(z_s^{t+1}) - g(z^*)),\\
            \langle z_s^{t+\nicefrac{1}{2}} - z_s^t + \gamma F_{\pi_s^t}(z_s^t), z_s^{t+1} - z_s^{t+\nicefrac{1}{2}}\rangle &\geqslant \gamma(g(z_s^{t+\nicefrac{1}{2}}) - g(z_s^{t+1})).
        \end{align*}
        Summing the inequalities, we obtain:
        \begin{align*}
        \gamma(g(z_s^{t+\nicefrac{1}{2}}) - g(z^*)) & \leqslant 
            \langle z_s^{t+1} - z_s^t, z^* - z_s^{t+1}\rangle+ \langle z_s^{t+\nicefrac{1}{2}} - z_s^t, z_s^{t+1} - z_s^{t+\nicefrac{1}{2}}\rangle  \\
            &\quad +\gamma\langle F_{\pi_s^t}(z_s^{t+\nicefrac{1}{2}}), z^* - z_s^{t+1}\rangle + \gamma\langle F_{\pi_s^t}(z_s^t), z_s^{t+1} - z_s^{t+\nicefrac{1}{2}}\rangle .
        \end{align*}
        Now, we add and subtract $z_s^{t+\nicefrac{1}{2}}$ to the right part of the third scalar product. Thus, rearranging terms, we arrive at
        \begin{align}
           \notag\gamma(g(z_s^{t+\nicefrac{1}{2}}) - g(z^*)) & \leqslant  \langle z_s^{t+1} - z_s^t, z^* - z_s^{t+1}\rangle + \langle z_s^{t+\nicefrac{1}{2}} - z_s^t, z_s^{t+1} - z_s^{t+\nicefrac{1}{2}}\rangle \\
           \notag& \quad+\gamma\langle F_{\pi_s^t}(z_s^{t+\nicefrac{1}{2}}), z^* - z_s^{t+\nicefrac{1}{2}}\rangle \\
           \label{t1:ineq5}& \quad+ \gamma\langle F_{\pi_s^t}(z_s^{t+\nicefrac{1}{2}}) - F_{\pi_s^t}(z_s^t), z_s^{t+\nicefrac{1}{2}} - z_s^{t+1}\rangle .
        \end{align}
        We want to rewrite the first two scalar products.  We use \eqref{bi:quadr}. Thus, we arrive at
        \begin{align*}
            2\langle z_s^{t+1} - z_s^t, z^* - z_s^{t+1}\rangle &= \|z_s^{t} - z^*\|^2 - \| z_s^{t+1} - z_s^t\|^2 - \|z^* - z_s^{t+1}\|^2,\\
            2\langle z_s^{t+\nicefrac{1}{2}} - z_s^t, z_s^{t+1} - z_s^{t+\nicefrac{1}{2}}\rangle &= \|z_s^{t+1} - z_s^t\|^2 - \|z_s^{t+\nicefrac{1}{2}} - z_s^t\|^2 - \|z_s^{t+1} - z_s^{t+\nicefrac{1}{2}}\|^2.
        \end{align*}
        Substituting this into \eqref{t1:ineq5}, we obtain:
        \begin{align*}
            \|z_s^{t+1} - z^*\|^2 &\leqslant \|z_s^t - z^*\|^2 - 2\gamma\left(\langle F_{\pi_s^t}(z_s^{t+\nicefrac{1}{2}}), z_s^{t+\nicefrac{1}{2}} - z^*\rangle + g(z_s^{t+\nicefrac{1}{2}}) - g(z^*) \right) \\
            & \quad + 2\gamma\langle F_{\pi_s^t}(z_s^{t+\nicefrac{1}{2}}) - F_{\pi_s^t}(z_s^t), z_s^{t+\nicefrac{1}{2}} - z_s^{t+1}\rangle \\
            & \quad - \|z_s^{t+\nicefrac{1}{2}} - z_s^t\|^2 - \|z_s^{t+1} - z_s^{t+\nicefrac{1}{2}}\|^2.
        \end{align*}
        Now, applying \eqref{bi:1} and Assumption \ref{as:lipschitz} to the second scalar product, we obtain:
        \begin{align}
            \notag\|z_s^{t+1} - z^*\|^2 &~~\leqslant ~\|z_s^t - z^*\|^2 - 2\gamma\bigl(\langle F_{\pi_s^t}(z_s^{t+\nicefrac{1}{2}}), z_s^{t+\nicefrac{1}{2}} - z^*\rangle\\
            \notag
            & ~~\quad + ~ g(z_s^{t+\nicefrac{1}{2}}) - g(z^*) \bigr) + \gamma^2\|F_{\pi_s^t}(z_s^{t+\nicefrac{1}{2}}) - F_{\pi_s^t}(z_s^t)\|^2\\
            \notag
            & ~~\quad + ~ \|z_s^{t+1} - z_s^{t+\nicefrac{1}{2}}\|^2 - \|z_s^{t+\nicefrac{1}{2}} - z_s^t\|^2 - \|z_s^{t+1} - z_s^{t+\nicefrac{1}{2}}\|^2\\
            \notag
            &\overset{\text{Ass.} \ref{as:lipschitz}}{\leqslant} \|z_s^t - z^*\|^2 + \left(\gamma^2 L^2 - 1\right) \|z_s^{t+\nicefrac{1}{2}} - z_s^t\|^2 \\
            \label{t1:ineq1}
            &~~~~~ - 2\gamma\left(\underbrace{\langle F_{\pi_s^t}(z_s^{t+\nicefrac{1}{2}}), z_s^{t+\nicefrac{1}{2}} - z^*\rangle + g(z_s^{t+\nicefrac{1}{2}}) - g(z^*)}_{T_1} \right).
        \end{align}
        To estimate the $T_1$ term, we compute the expectation:
        \begin{align*}
            \mathbb{E} T_1 & ~~~ = ~~~\mathbb{E}\langle F_{\pi_s^t}(z_s^{t+\nicefrac{1}{2}}), z_s^{t+\nicefrac{1}{2}} - z^*\rangle +  g(z_s^{t+\nicefrac{1}{2}}) -  g(z^*) \\
            &~~~= ~~~\mathbb{E} \langle F_{\pi_s^t}(z_s^{t+\nicefrac{1}{2}}) - F_{\pi_s^t}(z^*), z_s^{t+\nicefrac{1}{2}} - z^*\rangle \\
            & ~~~ \quad\quad + ~\mathbb{E}\langle F_{\pi_s^t}(z^*), z_s^{t+\nicefrac{1}{2}} - z^*\rangle + g(z_s^{t+\nicefrac{1}{2}}) -g(z^*)\\
            &~\overset{\text{Ass.} \ref{as:monotone}}{\geqslant} ~~\mu \mathbb{E} \|z_s^{t+\nicefrac{1}{2}} - z^*\|^2 + \mathbb{E} \langle F_{\pi_s^t}(z^*), z_s^{t+\nicefrac{1}{2}} - z_s^0\rangle \\
            & ~~~ \quad\quad + ~ \mathbb{E} \langle F_{\pi_s^t}(z^*), z_s^0 - z^*\rangle + g(z_s^{t+\nicefrac{1}{2}}) - g(z^*).
        \end{align*}
        Now we focus on the second scalar product. Utilizing the tower property and the unbiasedness of the stochastic operator at the points $z_s^0$ and $z^*$, we obtain $\mathbb E\left[\mathbb E_t\left[F_{\pi_s^t}(z^*)|z_s^0 - z^*\right]\right] = F(z^*)$. Thus, we continue the estimation of $\mathbb E T_1$:
        \begin{align*}
            \mathbb{E} T_1 &\overset{\eqref{bi:1}}{\geqslant} \mu \mathbb{E} \|z_s^{t+\nicefrac{1}{2}} - z^*\|^2 -  \frac{\gamma}{2\beta}\mathbb{E}\|F_{\pi_s^t}(z^*)\|^2 -  \frac{\beta}{2\gamma}\mathbb{E}\|z_s^{t+\nicefrac{1}{2}} - z_s^0\|^2 \\
            &~~ \quad ~~~+  \langle F(z^*), z_s^0 - z^*\rangle + g(z_s^{t+\nicefrac{1}{2}}) - g(z^*) \\
            &~~~ = \quad  \mu\mathbb{E}  \|z_s^{t+\nicefrac{1}{2}} - z^*\|^2 - \frac{\gamma}{2\beta}\mathbb{E}\|F_{\pi_s^t}(z^*)\|^2 - \frac{\beta}{2\gamma}\mathbb{E}\|z_s^{t+\nicefrac{1}{2}} - z_s^0\|^2\\
            &~~ \quad ~~~+ \langle F(z^*), z_s^0 - z_s^{t+\nicefrac{1}{2}}\rangle + \underbrace{\langle F(z^*), z_s^{t+\nicefrac{1}{2}} - z^*\rangle + g(z_s^{t+\nicefrac{1}{2}}) - g(z^*)}_{\geqslant 0 ~\eqref{eq:vi_setting}}\\
            &\overset{\eqref{bi:1}}{\geqslant} \mu \mathbb{E} \|z_s^{t+\nicefrac{1}{2}} - z^*\|^2 - \frac{\gamma}{2\beta} \mathbb{E}\|F_{\pi_s^t}(z^*)\|^2 \\
            & ~~ \quad ~~~- \frac{\gamma}{2\beta}\|F(z^*)\|^2 - \frac{\beta}{\gamma}\mathbb{E}\|z_s^{t+\nicefrac{1}{2}} - z_s^0\|^2.
        \end{align*}
        Here we introduce $\beta > 0$, which we will define later. Substituting this inequality into \eqref{t1:ineq1} yields:
        \begin{align}
            \notag \mathbb E\|z_s^{t+1} - z^*\|^2  &~~ \quad \leqslant \quad ~ \mathbb E\|z_s^t - z^*\|^2 - 2\gamma\mu \mathbb E\|z_s^{t+\frac{1}{2}} - z^*\|^2 \\
            \notag
            &\quad \quad \quad ~~+ (\gamma^2 L^2 - 1)\mathbb E\|z_s^{t+\frac{1}{2}} - z_s^t\|^2 + \frac{\gamma^2}{\beta}\mathbb E \|F_{\pi_s^t}(z^*)\|^2 \\
            \notag
            & \quad \quad \quad ~~+ \frac{\gamma^2}{\beta} \|F(z^*)\|^2 + 2\beta\mathbb E\|z_s^{t+\frac{1}{2}} - z_s^0\|^2\\
            \notag
            &\overset{(\text{Ass.} \ref{as:bound}, \ref{bi:quadrineq})}{\leqslant}(1 - \gamma\mu)\mathbb E\|z_s^t - z^*\|^2 + \frac{2\gamma^2}{\beta}\sigma_*^2 + 2\beta\mathbb E\|z_s^{t+\frac{1}{2}} - z_s^0\|^2 \\
            \label{t1:ineq2}
            &\quad \quad \quad ~~+ (\gamma^2 L^2 + 2\gamma\mu - 1)\mathbb E\|z_s^{t+\frac{1}{2}} - z_s^t\|^2.
        \end{align}
        Now, we evaluate the $\|z_s^{t+\frac{1}{2}} - z_s^0\|^2$ term:
        \allowdisplaybreaks
        \begin{align}
            \notag\|z_s^{t+\frac{1}{2}} - z_s^0\|^2 &~~\leqslant \left(1 + \frac{1}{a}\right)\|z_s^{t + \frac{1}{2}} - z_s^t\|^2 + (1 + a)\|z_s^t - z_s^0\|^2 \\
            \notag &~~\leqslant \left(1 + \frac{1}{a}\right)\|z_s^{t + \frac{1}{2}} - z_s^t\|^2 + (1 + a)\left(1 + \frac{1}{b}\right)\|z_s^t - z_s^{t - \frac{1}{2}}\|^2 \\
            \notag & ~~\quad + (1 + a)(1 + b)\|z_s^{t - \frac{1}{2}} - z_s^0\|^2 \\
            \notag &~~= \left(1 + \frac{1}{a}\right)\|z_s^{t + \frac{1}{2}} - z_s^t\|^2 \\
            \notag &~~\quad + (1 + a)\left(1 + \frac{1}{b}\right)\Bigl\|\text{prox}_{\gamma g}\left(z_s^{t-1} -\gamma F_{\pi_s^t}(z_s^{t - \frac{1}{2}})\right) \\
            \notag & ~~\quad- \text{prox}_{\gamma g}\left(z_s^{t-1} - \gamma F_{\pi_s^t}(z_s^{t - 1})\right)\Bigr\|^2 + (1 + a)(1 + b)\|z_s^{t - \frac{1}{2}} - z_s^0\|^2 \\
            \notag &\overset{\eqref{bi:Prox}}{\leqslant} \left(1 + \frac{1}{a}\right)\|z_s^{t + \frac{1}{2}} - z_s^t\|^2 \\
            \notag &~~\quad + (1 + a)\left(1 + \frac{1}{b}\right)\left\|\gamma F_{\pi_s^t}(z_s^{t - \frac{1}{2}}) - \gamma F_{\pi_s^t}(z_s^{t - 1})\right\|^2 \\
            \notag &~~\quad + (1 + a)(1 + b)\|z_s^{t - \frac{1}{2}} - z_s^0\|^2 \\
            \notag &~~\leqslant \left(1 + \frac{1}{a}\right)\|z_s^{t + \frac{1}{2}} - z_s^t\|^2 \\
            \notag &~~\quad + (1 + a)\left(1 + \frac{1}{b}\right)\gamma^2 L^2\|z_s^{t - \frac{1}{2}} - z_s^{t - 1}\|^2 \\
            \notag &~~\quad + (1 + a)(1 + b)\|z_s^{t - \frac{1}{2}} - z_s^0\|^2 \\
            \notag &~~\leqslant \left(1 + \frac{1}{a}\right)\|z_s^{t+\frac{1}{2}} - z_s^t\|^2 + \sum\limits_{i = 1}^{t-1}\|z_s^{i+\frac{1}{2}} - z_s^i\|^2  \\
            \notag &~~\quad ~ \cdot \left(\left(1+\frac{1}{a}\right)(1+a)(1+b) + (1+a)\left(1+\frac{1}{b}\right)\gamma^2 L^2 \right)\\
            \notag &~~\quad ~ \cdot \left[(1+a)(1+b)\right]^{t-1-i} \\
            \label{t1:ineq6}& ~~\quad + \left((1+a)\left(1+\frac{1}{b}\right)  + 1\right)\left[(1+a)(1+b)\right]^t\|z_s^{\frac{1}{2}} - z_s^0\|^2.
        \end{align}
        We choose $a = b = \frac{1}{n}$ and consider the coefficients for all three terms.
        \begin{align*}
            1 + \frac{1}{a} &= 1 + n,\\
            \left(\left(1 + \frac{1}{a}\right) + \frac{1}{b}\gamma^2 L^2\right)\left[(1+a)(1+b)\right]^{t-i} &= (1 + n + n\gamma^2 L^2)\left(1 + \frac{1}{n}\right)^{2(t-i)},\\
            \left((1+a)\left(1+\frac{1}{b}\right)+1\right)\left[(1+a)(1+b)\right]&^{t-i}\bigg|_{i = 0} \\
            &= \left(3 + \frac{1}{n} + n\right)\left(1 + \frac{1}{n}\right)^{2(t-i)}\bigg|_{i = 0}.
        \end{align*}
        We can evaluate the smaller terms from above by the largest one and consolidate them into a single sum. Thus, \eqref{t1:ineq6} transforms to
        \begin{equation*}
            \|z_s^{t+\frac{1}{2}} - z_s^0\|^2 \leqslant\sum\limits_{i = 0}^t \|z_s^{i + \frac{1}{2}} - z_s^i\|^2\left(1 + n + n\gamma^2 L^2\right)\left(1 + \frac{1}{n}\right)^{2(t-i)}.
        \end{equation*}
      Let us substitute the obtained inequality into \eqref{t1:ineq2}:
        \begin{align}
            \notag
            \mathbb{E} \|z_s^{t+1} - z^*\|^2 & \leqslant (1 - \gamma\mu)\mathbb{E} \|z_s^t - z^*\|^2 + (\gamma^2 L^2 + 2\gamma\mu - 1)\mathbb{E} \|z_s^{t+\frac{1}{2}} - z_s^t\|^2 ~ \\
            \label{t1:ineq4}
            & \quad + \frac{2\gamma^2}{\beta}\sigma_*^2 + 2\beta \sum\limits_{i = 0}^t \mathbb{E} \|z_s^{i + \frac{1}{2}} - z_s^i\|^2\left(1 + n + n\gamma^2 L^2\right)\left(1 + \frac{1}{n}\right)^{2(t-i)}.
        \end{align}
        Now we define a new sequence that contains iteration points in all epochs:
        \begin{equation*}
            \widetilde{z}_k = z_{t+sn}.
        \end{equation*}
        Thus, additionally considering $\left(1+\frac{1}{n}\right)^{2(t-i)} \leqslant \left(1+\frac{1}{n}\right)^{2n} \leqslant e^2 \leqslant 8$, we can rewrite \eqref{t1:ineq4} in the following form:
        \begin{align*}
            \mathbb{E} \|\widetilde{z}_{k+1} - z^*\|^2 & \leqslant (1 - \gamma\mu)\mathbb{E} \|\widetilde{z}_k - z^*\|^2 + (\gamma^2 L^2 + 2\gamma\mu - 1)\mathbb{E} \|\widetilde{z}_{k+\frac{1}{2}} - \widetilde{z}_k\|^2 ~ \\
            &\quad + \frac{2\gamma^2}{\beta}\sigma_*^2 + 16\beta \sum\limits_{i = 0}^n \mathbb{E} \|\widetilde{z}_{k - i + \frac{1}{2}} - \widetilde{z}_{k-i}\|^2\left(1 + n + n\gamma^2 L^2\right).
        \end{align*}
        Let us pay attention to the $\sum\limits_{i = 0}^n \mathbb{E} \|\widetilde{z}_{k - i + \frac{1}{2}} - \widetilde{z}_{k-i}\|^2$ term in the obtained inequality. For the original sequence, this term represented the sum of the norms from the beginning of the current epoch to the current iteration $t$ and could contain a maximum of $n$ terms. Thus, a new expression that includes the sum of $n$ norms up to the current iteration $k$ serves as an upper bound, confirming that our expression is correct.
        Now we define $p_k = p^k = \left(1 - \frac{\gamma\mu}{2}\right)^{-k}$ and summarize both sides over all iterations with coefficients $p_k$:
        \begin{align}
        \notag
        \sum\limits_{k = 0}^{Sn-1} p_k\mathbb{E} \|\widetilde{z}_{k+1} - z^*\|^2 & \leqslant (1 - \gamma\mu)\sum\limits_{k = 0}^{Sn-1} p_k\mathbb{E} \|\widetilde{z}_k - z^*\|^2 \\
        \notag
        & \quad + (\gamma^2 L^2 + 2\gamma\mu - 1)\sum\limits_{k = 0}^{Sn-1} p_k\mathbb{E} \|\widetilde{z}_{k+\frac{1}{2}} - \widetilde{z}_k\|^2 \\
        \notag
        & \quad + \sum\limits_{k = 0}^{Sn-1} p_k\sum\limits_{i = 0}^n \mathbb{E} \|\widetilde{z}_{k - i + \frac{1}{2}} - \widetilde{z}_{k-i}\|^2\\
        \label{t1:ineq3}
        & \quad \cdot 16\beta\left(1 + n + n\gamma^2 L^2\right) + \frac{2\gamma^2\sigma_*^2}{\beta}\sum\limits_{k = 0}^{Sn-1}p_k.
        \end{align}
        Now we need to estimate the following term:
        \begin{align*}
            \sum\limits_{k = 0}^{Sn-1} p_k\sum\limits_{i = 0}^n \mathbb{E} \|\widetilde{z}_{k - i + \frac{1}{2}} - \widetilde{z}_{k-i}\|^2 &\leqslant p_n \sum\limits_{k = 0}^{Sn-1} \sum\limits_{i = 0}^n p_{k-i}\mathbb{E} \|\widetilde{z}_{k - i + \frac{1}{2}} - \widetilde{z}_{k-i}\|^2 \\
            &\leqslant p_n n \sum\limits_{k = 0}^{Sn-1} p_k\mathbb{E} \|\widetilde{z}_{k + \frac{1}{2}} - \widetilde{z}_k\|^2.
        \end{align*}
        Note that we define points $\widetilde{z}_{-n}, \widetilde{z}_{-n+\frac{1}{2}}, \ldots, \widetilde{z}_{-\frac{1}{2}}$ by shifting the sequence $\{\widetilde{z}_k\}$ on $n$ points. Since $p_n = \left(1 - \frac{\gamma\mu}{2}\right)^{-n} = \left(1 - \frac{\gamma\mu n}{2n}\right)^{-n} \leqslant e^{\frac{\gamma\mu n}{2}}$, we choose $\gamma \leqslant \frac{1}{2\mu n}$ and obtain $p_n \leqslant e^{\frac{1}{4}} \leqslant 2$. Substituting this into \eqref{t1:ineq3}, we obtain:
        \begin{align*}
            \sum\limits_{k = 0}^{Sn-1} p_k\mathbb{E} \|\widetilde{z}_{k+1} - z^*\|^2 & \leqslant (1 - \gamma\mu)\sum\limits_{k = 0}^{Sn-1} p_k\mathbb{E} \|\widetilde{z}_k - z^*\|^2 + \frac{2\gamma^2\sigma_*^2}{\beta}\sum\limits_{k = 0}^{Sn-1}p_k \\
            & \quad + (\gamma^2 L^2 + 2\gamma\mu - 1)\sum\limits_{k = 0}^{Sn-1} p_k\mathbb{E} \|\widetilde{z}_{k+\frac{1}{2}} - \widetilde{z}_k\|^2 \\
            & \quad + 32\beta\left(1 + n + n\gamma^2 L^2\right)n\sum\limits_{k = 0}^{Sn-1} p_k \mathbb{E} \|\widetilde{z}_{k + \frac{1}{2}} - \widetilde{z}_k\|^2.
        \end{align*}
        We consider the coefficient before $\sum\limits_{k = 0}^{Sn-1} p_k\mathbb E\|\widetilde{z}_{k+\frac{1}{2}} - \widetilde{z}_k\|^2$ and we make it negative  by selecting $\gamma$ and $\beta$.
        \begin{gather*}
            32\beta(1 + n + n\gamma^2 L^2)n + \gamma^2 L^2 + 2\gamma\mu - 1 \leqslant 0
        \end{gather*}
        We need $\gamma\leqslant\frac{1}{6L}, \beta = \frac{1}{64n^2}.$ Then to satisfy the previous estimate on gamma we finally put $\gamma \leqslant \min\left\{\frac{1}{2\mu n}, \frac{1}{6L}\right\}$ and, assuming $n > 3$, have $$\frac{1}{2n} + \frac{1}{2} + \frac{1}{72} + \frac{1}{36} + \frac{1}{3} - 1 \leqslant 0.$$
        In this way,
        \begin{equation*}
            \sum\limits_{k = 0}^{Sn-1} p_k\mathbb{E} \|\widetilde{z}_{k+1} - z^*\|^2\leqslant (1 - \gamma\mu)\sum\limits_{k = 0}^{Sn-1} p_k\mathbb{E} \|\widetilde{z}_k - z^*\|^2 +  \frac{2\gamma^2\sigma_*^2}{\beta}\sum\limits_{k = 0}^{Sn-1}p_k.
        \end{equation*}
        Thus, substituting definition of $p_t$, we obtain:
        \allowdisplaybreaks
        \begin{align*}
            \sum\limits_{k = 0}^{Sn-1} \left(1-\frac{\gamma\mu}{2}\right)^{-k}\mathbb{E} \|\widetilde{z}_{k+1} - z^*\|^2 & \leqslant \sum\limits_{k = 0}^{Sn-1} \left(1-\frac{\gamma\mu}{2}\right)^{-k + 1}\mathbb{E} \|\widetilde{z}_k - z^*\|^2 \\
            & \quad + \frac{2\gamma^2\sigma_*^2}{\beta}\sum\limits_{k = 0}^{Sn-1}\left(1-\frac{\gamma\mu}{2}\right)^{-k},\\
            \left(1-\frac{\gamma\mu}{2}\right)^{-(Sn-1)}\mathbb{E} \|\widetilde{z}_{Sn} - z^*\|^2 & \leqslant \left(1-\frac{\gamma\mu}{2}\right)\mathbb{E} \|\widetilde{z}_0 - z^*\|^2 \\
            & \quad + \frac{2\gamma^2\sigma_*^2}{\beta}\sum\limits_{k = 0}^{Sn-1}\left(1-\frac{\gamma\mu}{2}\right)^{-k},\\
            \mathbb{E} \|\widetilde{z}_{Sn} - z^*\|^2 & \leqslant \left(1-\frac{\gamma\mu}{2}\right)^{Sn}\mathbb{E} \|\widetilde{z}_0 - z^*\|^2\\
            & \quad + \frac{2\gamma^2\sigma_*^2}{\beta}\sum\limits_{k = 0}^{Sn-1}\left(1-\frac{\gamma\mu}{2}\right)^{Sn-k-1}\\
            & = \left(1-\frac{\gamma\mu}{2}\right)^{Sn}\mathbb{E} \|\widetilde{z}_0 - z^*\|^2\\
            & \quad + \frac{2\gamma^2\sigma_*^2}{\beta}\sum\limits_{k = 0}^{Sn-1}\left(1-\frac{\gamma\mu}{2}\right)^{k}.
        \end{align*}
        Finally, estimating the geometric progression in the last term as $\sum\limits_{k = 0}^{Sn-1}\left(1-\frac{\gamma\mu}{2}\right)^{k} \leqslant \frac{2}{\gamma\mu}$, we can write the final statement of the theorem:
        \begin{align*}
            \mathbb{E} \|z_S^n - z^*\|^2 & \leqslant \left(1-\frac{\gamma\mu}{2}\right)^{Sn}\mathbb{E} \|z_0^0 - z^*\|^2 + \frac{4\gamma\sigma_*^2}{\beta\mu}\\
            & = \left(1-\frac{\gamma\mu}{2}\right)^{Sn}\mathbb{E} \|z_0^0 - z^*\|^2 + \frac{256\gamma n^2\sigma_*^2}{\mu}.
        \end{align*}
    \end{proof}
\textbf{Corollary \ref{cor:eg}.}
\textit{Suppose Assumptions \ref{as:lipschitz}, \ref{as:monotone} hold. Then Algorithms \ref{alg:rrextragrad}, \ref{alg:soextragrad} with \\ $\gamma \leqslant \min\left\{\frac{1}{2\mu n}, \frac{1}{6L}, \frac{2\log\left(\max\left\{2, \frac{\mu^2\|z_0^0 - z^*\|^2 T}{512 n^2\sigma^2_*}\right\}\right)}{\mu T}\right\}$, to reach $\varepsilon$-accuracy, where $\varepsilon \sim \|z_S^n - z^*\|^2$, needs}
    \begin{equation*}
    \mathcal{\widetilde{O}}\left(\left(n + \frac{L}{\mu}\right)\log\left(\frac{1}{\varepsilon}\right) + \frac{n^2\sigma^2_*}{\mu^2\varepsilon}\right) ~~\textit{iterations and oracle calls.}
    \end{equation*}
\begin{proof}
    For the result obtained in Theorem \ref{th:eg}, we utilize Lemma 2 from \cite{stich2019unified} and, using special tuning of $\gamma$, such as $\gamma \leqslant \min\left\{\frac{1}{2\mu n}, \frac{1}{6L}, \frac{2\log\left(\max\left\{2, \frac{\mu^2\|z_0^0 - z^*\|^2 T}{512 n^2\sigma^2_*}\right\}\right)}{\mu T}\right\}$, we obtain that we need $\mathcal{\widetilde{O}}\left(\left(n + \frac{L}{\mu}\right)\log\left(\frac{1}{\varepsilon}\right) + \frac{n^2\sigma^2_*}{\mu^2\varepsilon}\right)$ iterations and oracle calls to reach $\varepsilon$-accuracy, where $\varepsilon \sim \|z_S^n - z^*\|^2$.
\end{proof}

\section{Extragradient with Variance Reduction}\label{D}
\textbf{Theorem \ref{th:proxegvr}.}
\textit{Suppose that Assumptions \ref{as:lipschitz}, \ref{as:monotone} hold. Then for Algorithm \ref{alg:proxextragradvr} with $\gamma \leqslant\frac{(1-\alpha)\mu}{6L^2}, p = \frac{1}{n}$ and $V_s^t =\mathbb E \|z_s^t - z^*\|^2 + \mathbb E\|\omega_s^t - z^*\|^2$, after $T$ iterations we have}
    \begin{equation*}
        V_S^n \leqslant \left(1 - \frac{\gamma\mu}{4}\right)^T V_0^0.
    \end{equation*}
\begin{proof}
    We start with substituting both steps of Algorithm \ref{alg:proxextragradvr} to \eqref{ineq:prox}:
        \begin{align*}
            \langle z_s^{t+1} - \overline{z}_s^t + \gamma \hat{F}(z_s^{t+\nicefrac{1}{2}}), z^* - z_s^{t+1}\rangle &\geqslant \gamma (g(z_s^{t+1}) - g(z^*)),\\
            \langle z_s^{t+\nicefrac{1}{2}} - \overline{z}_s^t + \gamma F(\omega_s^t), z_s^{t+1} - z_s^{t+\nicefrac{1}{2}}\rangle &\geqslant \gamma (g(z_s^{t+\nicefrac{1}{2}}) - g(z_s^{t+1})).
        \end{align*}
        Let us summarize this two inequalities:
        \begin{align*}
        \gamma (g(z_s^{t+\nicefrac{1}{2}}) - g(z^*)) &\leqslant
            \langle z_s^{t+1} - \overline{z}_s^t, z^* - z_s^{t+1}\rangle+ \langle z_s^{t+\nicefrac{1}{2}} - \overline{z}_s^t, z_s^{t+1} - z_s^{t+\nicefrac{1}{2}}\rangle\\
            & \quad + \gamma\langle \hat{F}(z_s^{t+\nicefrac{1}{2}}), z^* - z_s^{t+1}\rangle + \gamma\langle F(\omega_s^t), z_s^{t+1} - z_s^{t+\nicefrac{1}{2}}\rangle .
        \end{align*}
        Now, we add and subtract $z_s^{t+\nicefrac{1}{2}}$ to the right part of the third scalar product. Thus, rearranging terms and utilizing the definition of $\hat{F}(z_s^{t+\nicefrac{1}{2}})$, we arrive at:
        \begin{align}
            \notag
            &\underbrace{\langle z_s^{t+1} - \overline{z}_s^t, z^* - z_s^{t+1}\rangle}_{T_1} + \underbrace{\langle z_s^{t+\nicefrac{1}{2}} - \overline{z}_s^t, z_s^{t+1} - z_s^{t+\nicefrac{1}{2}}\rangle}_{T_2}\\
            \notag
            &\quad + \underbrace{\gamma\langle F_{\pi_s^t}(\omega_s^t) - F_{\pi_s^t}(z_s^{t+\nicefrac{1}{2}}), z_s^{t+1} - z_s^{t+\nicefrac{1}{2}}\rangle}_{T_3} \\
            \label{t2:ineq1}
            &\quad + \underbrace{\gamma\langle \hat{F}(z_s^{t+\nicefrac{1}{2}}), z^* - z_s^{t+\nicefrac{1}{2}}\rangle + \gamma (g(z^*) - g(z_s^{t+\nicefrac{1}{2}}))}_{T_4} \geqslant 0.
        \end{align}
        We defined terms as $T_1, T_2, T_3, T_4$, respectively. Let us estimate them separately. We start with $T_1$ and $T_2$. To estimate them, firstly, we use the definition of $\overline{z}_s^t$ and, secondly, use \eqref{bi:quadr}. Thus, we obtain:
        \begin{align*}
            2T_1 &= 2\langle z_s^{t+1} - \overline{z}_s^t, z^* - z_s^{t+1}\rangle \\
            &= 2\alpha \langle z_s^{t+1} - z_s^t, z^* - z_s^{t+1}\rangle + 2(1-\alpha)\langle z_s^{t+1} - \omega_s^t, z^* - z_s^{t+1}\rangle\\
            &= \alpha (\|z^* - z_s^t\|^2 - \|z_s^{t+1} - z_s^t\|^2 - \|z^* - z_s^{t+1}\|^2) \\
            &\quad + (1-\alpha)(\|z^* - \omega_s^t\|^2 - \|z_s^{t+1} - \omega_s^t\|^2 - \|z^* - z_s^{t+1}\|^2)\\
            &= \alpha \|z_s^t - z^*\|^2 - \|z_s^{t+1} - z^*\|^2 + (1-\alpha)\|\omega_s^t - z^*\|^2 \\
            &\quad - \alpha \|z_s^{t+1} - z_s^t\|^2 - (1-\alpha) \|z_s^{t+1} - \omega_s^t\|^2.
        \end{align*}
        The same holds for $T_2$:
        \begin{align*}
            2T_2 &= 2\langle z_s^{t+\nicefrac{1}{2}} - \overline{z}_s^t, z_s^{t+1} - z_s^{t+\nicefrac{1}{2}}\rangle \\
            &= 2\alpha \langle z_s^{t+\nicefrac{1}{2}} - z_s^t, z_s^{t+1} - z_s^{t+\nicefrac{1}{2}}\rangle + 2(1-\alpha)\langle z_s^{t+\nicefrac{1}{2}} - \omega_s^t, z_s^{t+1} - z_s^{t+\nicefrac{1}{2}}\rangle\\
            &= \alpha (\|z_s^{t+1} - z_s^t\|^2 - \|z_s^{t+\nicefrac{1}{2}} - z_s^t\|^2 - \|z_s^{t+1} - z_s^{t+\nicefrac{1}{2}}\|^2) \\
            &\quad + (1-\alpha)(\|z_s^{t+1} - \omega_s^t\|^2 - \|z_s^{t+\nicefrac{1}{2}} - \omega_s^t\|^2 - \|z_s^{t+1} - z_s^{t+\nicefrac{1}{2}}\|^2)\\
            &= \alpha \|z_s^{t+1} - z_s^t\|^2 - \|z_s^{t+1} - z_s^{t+\nicefrac{1}{2}}\|^2 + (1-\alpha)\|z_s^{t+1} - \omega_s^t\|^2 \\
            &\quad - \alpha \|z_s^{t+\nicefrac{1}{2}} - z_s^t\|^2 - (1-\alpha) \|z_s^{t+\nicefrac{1}{2}} - \omega_s^t\|^2.
        \end{align*}
        Now, we moving to the estimate of $T_3$:
        \begin{align*}
            2T_3 & ~~~= ~~~2\gamma\langle F_{\pi_s^t}(\omega_s^t) - F_{\pi_s^t}(z_s^{t+\nicefrac{1}{2}}), z_s^{t+1} - z_s^{t+\nicefrac{1}{2}}\rangle \\
            &\overset{\eqref{bi:1}}{\leqslant} \frac{\gamma^2}{\tau}\|F_{\pi_s^t}(\omega_s^t) - F_{\pi_s^t}(z_s^{t+\nicefrac{1}{2}})\|^2 + \tau \|z_s^{t+1} - z_s^{t+\nicefrac{1}{2}}\|^2 \\
            &~\overset{\eqref{bi:CauchySchwarz}}{\leqslant} ~\frac{\gamma^2 L^2}{\tau}\|z_s^{t+\nicefrac{1}{2}} - \omega_s^t\|^2 + \tau \|z_s^{t+1} - z_s^{t+\nicefrac{1}{2}}\|^2;
        \end{align*}
        Here we introduced $\tau > 0$, which we will define later. Last, we do the same for $T_4$:
        \begin{align*}
            2T_4 &~\quad \quad=\quad \quad  2\gamma\langle \hat{F} (z_s^{t+\nicefrac{1}{2}}), z^* - z_s^{t+\nicefrac{1}{2}}\rangle + 2\gamma (g(z^*) - g(z_s^{t+\nicefrac{1}{2}}))\\
            &~\quad \quad =\quad \quad  2\gamma\langle \hat{F} (z_s^{t+\nicefrac{1}{2}}) - F(z_s^{t+\nicefrac{1}{2}}), z^* - z_s^{t+\nicefrac{1}{2}}\rangle \\
            &\quad\quad \quad\quad\quad + 2\gamma\langle F (z_s^{t+\nicefrac{1}{2}}) - F(z^*), z^* - z_s^{t+\nicefrac{1}{2}}\rangle \\
            &\quad \quad \quad \quad \quad + 2\gamma\left(\underbrace{\langle F(z^*), z^* - z_s^{t+\nicefrac{1}{2}}\rangle + g(z^*) - g(z_s^{t+\nicefrac{1}{2}})}_{\leqslant 0 ~\eqref{eq:vi_setting}}\right) \\
            &\overset{(\ref{bi:1}, \text{Ass.} \ref{as:lipschitz})}{\leqslant} \frac{4\gamma^2 L^2}{\tau} \|z_s^{t+\nicefrac{1}{2}} - \omega_s^t\|^2 + \tau \|z_s^{t+\nicefrac{1}{2}} - z^*\|^2 \\
            &\quad \quad \quad \quad \quad  - 2\gamma\langle F (z_s^{t+\nicefrac{1}{2}}) - F(z^*), z_s^{t+\nicefrac{1}{2}} - z^*\rangle\\
            &~\quad ~ \overset{(\text{Ass.} \ref{as:monotone})}{\leqslant} \quad ~ \frac{4\gamma^2 L^2}{\tau} \|z_s^{t+\nicefrac{1}{2}} - \omega_s^t\|^2 + \tau \|z_s^{t+\nicefrac{1}{2}} - z^*\|^2 - 2\gamma\mu\|z_s^{t+\nicefrac{1}{2}} - z^*\|^2. 
        \end{align*}
       Substituting all the obtained estimates into \eqref{t2:ineq1}, we arrive at
        \begin{align*}
            0 &\leqslant \alpha \|z_s^t - z^*\|^2 - \|z_s^{t+1} - z^*\|^2 + (1-\alpha)\|\omega_s^t - z^*\|^2 - \alpha \|z_s^{t+1} - z_s^t\|^2 \\
            &\quad - (1-\alpha) \|z_s^{t+1} - \omega_s^t\|^2 + \alpha \|z_s^{t+1} - z_s^t\|^2 - \|z_s^{t+1} - z_s^{t+\nicefrac{1}{2}}\|^2 \\
            & \quad + (1-\alpha)\|z_s^{t+1} - \omega_s^t\|^2 - \alpha \|z_s^{t+\nicefrac{1}{2}} - z_s^t\|^2 - (1-\alpha) \|z_s^{t+\nicefrac{1}{2}} - \omega_s^t\|^2 \\
            & \quad + \frac{\gamma^2 L^2}{\tau}\|z_s^{t+\nicefrac{1}{2}} - \omega_s^t\|^2 + \tau \|z_s^{t+1} - z_s^{t+\nicefrac{1}{2}}\|^2 +\frac{4\gamma^2 L^2}{\tau} \|z_s^{t+\nicefrac{1}{2}} - \omega_s^t\|^2 \\
            & \quad + \tau \|z_s^{t+\nicefrac{1}{2}} - z^*\|^2 - 2\gamma\mu\|z_s^{t+\nicefrac{1}{2}} - z^*\|^2.
        \end{align*}
        By grouping the coefficients of the same terms, we get:
        \begin{align}
        \notag
            \|z_s^{t+1} - z^*\|^2 &\leqslant \alpha \|z_s^t - z^*\|^2 + (1-\alpha)\|\omega_s^t - z^*\|^2 \\
            \notag
            & \quad + \left(\frac{5\gamma^2 L^2}{\tau} - (1-\alpha)\right) \|z_s^{t+\nicefrac{1}{2}} - \omega_s^t\|^2 - (1-\tau)\|z_s^{t+1} - z_s^{t+\nicefrac{1}{2}}\|^2 \\
            \label{t2:ineq2}
            & \quad - (2\gamma\mu - \tau)\|z_s^{t+\nicefrac{1}{2}} - z^*\|^2 - \alpha \|z_s^{t+\nicefrac{1}{2}} - z_s^t\|^2.
        \end{align}
        Now, we want to estimate the $- (2\gamma\mu - \tau)\|z_s^{t+\nicefrac{1}{2}} - z^*\|^2$ term. To do this we split it into two equal parts. To the first part we add and subtract $\omega_s^t$, and to the second -- $z_s^t$. After that we use \eqref{bi:quadrineq} for both terms:
        \begin{align*}
            - \left(2\gamma\mu - \tau\right)\|z_s^{t+\nicefrac{1}{2}} - z^*\|^2 &= - \left(\gamma\mu - \frac{\tau}{2}\right)\|z_s^{t+\nicefrac{1}{2}} - z^*\|^2 - \left(\gamma\mu - \frac{\tau}{2}\right)\|z_s^{t+\nicefrac{1}{2}} - z^*\|^2 \\
            &= - \left(\gamma\mu - \frac{\tau}{2}\right)\|z_s^{t+\nicefrac{1}{2}} - \omega_s^t + \omega_s^t - z^*\|^2 \\
            & \quad - \left(\gamma\mu - \frac{\tau}{2}\right)\|z_s^{t+\nicefrac{1}{2}} - z_s^t + z_s^t - z^*\|^2 \\
            & \leqslant \left(\gamma\mu - \frac{\tau}{2}\right) \|z_s^{t+\nicefrac{1}{2}} - \omega_s^t\|^2 - \left(\frac{\gamma\mu}{2} - \frac{\tau}{4}\right)\|\omega_s^t - z^*\|^2 \\
            &\quad+ \left(\gamma\mu - \frac{\tau}{2}\right) \|z_s^{t+\nicefrac{1}{2}} - z_s^t\|^2 - \left(\frac{\gamma\mu}{2} - \frac{\tau}{4}\right)\|z_s^t - z^*\|^2.
        \end{align*}
        Substituting this into \eqref{t2:ineq2},
        \begin{align*}
            \|z_s^{t+1} - z^*\|^2 &\leqslant \left(\alpha - \frac{\gamma\mu}{2} + \frac{\tau}{4}\right) \|z_s^t - z^*\|^2 + \left(1-\alpha - \frac{\gamma\mu}{2} + \frac{\tau}{4} \right)\|\omega_s^t - z^*\|^2 \\
            &\quad + \left(\frac{5\gamma^2 L^2}{\tau} + \gamma\mu - \frac{\tau}{2} - (1-\alpha)\right) \|z_s^{t+\nicefrac{1}{2}} - \omega_s^t\|^2 \\
            &\quad - (1-\tau)\|z_s^{t+1} - z_s^{t+\nicefrac{1}{2}}\|^2 + \left(\gamma\mu - \frac{\tau}{2} -\alpha \right)\|z_s^{t+\nicefrac{1}{2}} - z_s^t\|^2.
        \end{align*}
        We want to choose parameters such that coefficients before the last three terms would be non-positive. Let us start with $\|z_s^{t+\nicefrac{1}{2}} - \omega_s^t\|^2$ term. 
        \begin{align*}
            \text{We pick~}\tau &= \gamma\mu;\\
            \text{We want~}1- \alpha &\geqslant \frac{5\gamma L^2}{\mu} + \frac{\gamma \mu}{2};\\
            \text{It is enough for us that~}\gamma &\leqslant \frac{(1-\alpha)\mu}{6L^2}.
         \end{align*}
        Obviously, with this choice of $\gamma$ and $\alpha$, the last two terms are less than zero. In that way, we obtain:
        \begin{equation*}
            \|z_s^{t+1} - z^*\|^2 \leqslant \left(\alpha - \frac{\gamma\mu}{4}\right)\|z_s^t - z^*\|^2 + \left(1 - \alpha - \frac{\gamma\mu}{4}\right)\|\omega_s^t - z^*\|^2.
        \end{equation*}
        According to the condition for updating the point $\omega_s^t$,
        \begin{equation*}
            \mathbb E\|\omega_s^{t+1} - z^*\|^2 = p\|z_s^t - z^*\|^2 + (1-p)\|\omega_s^t - z^*\|^2.
        \end{equation*}
        In that way:
        \begin{align*}
            \mathbb E\|z_s^{t+1} - z^*\|^2 + \frac{1-\alpha}{p} \mathbb E\|\omega_s^{t+1} - z^*\|^2 &\leqslant \left(1 - \frac{\gamma\mu}{4}\right)\mathbb E\|z_s^t - z^*\|^2 \\
            & \quad + \left((1-\alpha)\left(1 + \frac{1}{p} - 1\right) - \frac{\gamma\mu}{4}\right)\mathbb E\|\omega_s^{t+1} - z^*\|^2.
        \end{align*}
        Now we put $\alpha = 1 - p$ and obtain:
        \begin{align*}
            \mathbb E\|z_s^{t+1} - z^*\|^2 + \mathbb E\|\omega_s^{t+1} - z^*\|^2 \leqslant \left(1 - \frac{\gamma\mu}{4}\right)\left(\mathbb E\|z_s^t - z^*\|^2 + \mathbb E\|\omega_s^t - z^*\|^2\right).
        \end{align*}
        Denoting $V_s^t = \mathbb E\|z_s^t - z^*\|^2 + \mathbb E\|\omega_s^t - z^*\|^2$ and going into recursion over all epochs and iterations, we get:
        \begin{equation*}
            V_S^n \leqslant \left(1 - \frac{\gamma\mu}{4}\right)^T V_0^0,
        \end{equation*}
        where $T$ is the total number of iterations. 
\end{proof}    
\textbf{Corollary \ref{cor:proxegvr}.}
\textit{Suppose that Assumptions \ref{as:lipschitz}, \ref{as:monotone} hold. Then Algorithm \ref{alg:proxextragradvr} with $\gamma \leqslant\frac{(1-\alpha)\mu}{6L^2}, p = \frac{1}{n}$ and $V_s^t =\mathbb E \|z_s^t - z^*\|^2 + \mathbb E\|\omega_s^t - z^*\|^2$, to reach $\varepsilon$-accuracy, where $\varepsilon \sim V_S^n$, needs}
\begin{align*}
    &\mathcal{O}\left(n\frac{L^2}{\mu^2}\log\left(\frac{1}{\varepsilon}\right)\right) ~~\text{iterations and oracle calls.}
    \end{align*}
\begin{proof}   
        Substituting estimation of $\gamma$ to the result of Theorem \ref{th:proxegvr} we obtain, that method to converge to $\varepsilon$-accuracy, where $\varepsilon = V_S^n$, needs $\mathcal{O}\left(\frac{L^2}{p\mu^2}\log\left(\frac{1}{\varepsilon}\right)\right)$ iterations. At the same time each iteration costs $pn + 2$ calls to $F_{\pi}$. Thus, we obtain $\mathcal{O}\left(\left(n\frac{L^2}{\mu^2} + \frac{L^2}{p\mu^2}\right)\log\left(\frac{1}{\varepsilon}\right)\right)$ oracle complexity. Finally, the optimal choice $p = \frac{1}{n}$ gives $\mathcal{O}\left(n\frac{L^2}{\mu^2}\log\left(\frac{1}{\varepsilon}\right)\right)$ iteration and oracle complexity. This ends the proof.
\end{proof}
    
\end{appendixpart}


\begin{thebibliography}{72}
\providecommand{\natexlab}[1]{#1}
\providecommand{\url}[1]{\texttt{#1}}
\expandafter\ifx\csname urlstyle\endcsname\relax
  \providecommand{\doi}[1]{doi: #1}\else
  \providecommand{\doi}{doi: \begingroup \urlstyle{rm}\Url}\fi

\bibitem[Alacaoglu and Malitsky(2022)]{alacaoglu2022stochastic}
Ahmet Alacaoglu and Yura Malitsky.
\newblock Stochastic variance reduction for variational inequality methods.
\newblock In \emph{Conference on Learning Theory}, pages 778--816. PMLR, 2022.

\bibitem[Alacaoglu et~al.(2021)Alacaoglu, Malitsky, and Cevher]{alacaoglu2021forward}
Ahmet Alacaoglu, Yura Malitsky, and Volkan Cevher.
\newblock Forward-reflected-backward method with variance reduction.
\newblock \emph{Computational Optimization and Applications}, 80\penalty0 (2):\penalty0 321--346, 2021.

\bibitem[Bach et~al.(2008)Bach, Mairal, and Ponce]{bach2008convex}
Francis Bach, Julien Mairal, and Jean Ponce.
\newblock Convex sparse matrix factorizations.
\newblock \emph{arXiv preprint arXiv:0812.1869}, 2008.

\bibitem[Ben-Tal et~al.(2009)Ben-Tal, El~Ghaoui, and Nemirovski]{ben2009robust}
Aharon Ben-Tal, Laurent El~Ghaoui, and Arkadi Nemirovski.
\newblock \emph{Robust optimization}, volume~28.
\newblock Princeton University Press, 2009.

\bibitem[Bertsekas and Tsitsiklis(2000)]{bertsekas2000gradient}
Dimitri~P Bertsekas and John~N Tsitsiklis.
\newblock Gradient convergence in gradient methods with errors.
\newblock \emph{SIAM Journal on Optimization}, 10\penalty0 (3):\penalty0 627--642, 2000.

\bibitem[Beznosikov et~al.(2020)Beznosikov, Samokhin, and Gasnikov]{beznosikov2020distributed}
Aleksandr Beznosikov, Valentin Samokhin, and Alexander Gasnikov.
\newblock Distributed saddle-point problems: Lower bounds, near-optimal and robust algorithms.
\newblock \emph{arXiv preprint arXiv:2010.13112}, 2020.

\bibitem[Beznosikov et~al.(2023)Beznosikov, Polyak, Gorbunov, Kovalev, and Gasnikov]{beznosikov2023smooth}
Aleksandr Beznosikov, Boris Polyak, Eduard Gorbunov, Dmitry Kovalev, and Alexander Gasnikov.
\newblock Smooth monotone stochastic variational inequalities and saddle point problems: A survey.
\newblock \emph{European Mathematical Society Magazine}, \penalty0 (127):\penalty0 15--28, 2023.

\bibitem[Bottou(2009)]{bottou2009curiously}
L{\'e}on Bottou.
\newblock Curiously fast convergence of some stochastic gradient descent algorithms.
\newblock In \emph{Proceedings of the Symposium on Learning and Data Science, Paris}, volume~8, pages 2624--2633. Citeseer, 2009.

\bibitem[Bottou(2010)]{bottou2010large}
L{\'e}on Bottou.
\newblock Large-scale machine learning with stochastic gradient descent.
\newblock In \emph{Proceedings of COMPSTAT'2010: 19th International Conference on Computational StatisticsParis France, August 22-27, 2010 Keynote, Invited and Contributed Papers}, pages 177--186. Springer, 2010.

\bibitem[Browder(1965)]{browder1965nonexpansive}
Felix~E Browder.
\newblock Nonexpansive nonlinear operators in a banach space.
\newblock \emph{Proceedings of the National Academy of Sciences}, 54\penalty0 (4):\penalty0 1041--1044, 1965.

\bibitem[Carmon et~al.(2019)Carmon, Jin, Sidford, and Tian]{carmon2019variance}
Yair Carmon, Yujia Jin, Aaron Sidford, and Kevin Tian.
\newblock Variance reduction for matrix games.
\newblock \emph{Advances in Neural Information Processing Systems}, 32, 2019.

\bibitem[Chambolle and Pock(2011)]{chambolle2011first}
Antonin Chambolle and Thomas Pock.
\newblock A first-order primal-dual algorithm for convex problems with applications to imaging.
\newblock \emph{Journal of Mathematical Imaging and Vision}, 40:\penalty0 120--145, 2011.

\bibitem[Chang and Lin(2011)]{chang2011libsvm}
Chih-Chung Chang and Chih-Jen Lin.
\newblock Libsvm: a library for support vector machines.
\newblock \emph{ACM Transactions on Intelligent Systems and Technology (TIST)}, 2\penalty0 (3):\penalty0 1--27, 2011.

\bibitem[Chavdarova et~al.(2019)Chavdarova, Gidel, Fleuret, and Lacoste-Julien]{chavdarova2019reducing}
Tatjana Chavdarova, Gauthier Gidel, Fran{\c{c}}ois Fleuret, and Simon Lacoste-Julien.
\newblock Reducing noise in gan training with variance reduced extragradient.
\newblock \emph{Advances in Neural Information Processing Systems}, 32, 2019.

\bibitem[Daskalakis et~al.(2017)Daskalakis, Ilyas, Syrgkanis, and Zeng]{daskalakis2017training}
Constantinos Daskalakis, Andrew Ilyas, Vasilis Syrgkanis, and Haoyang Zeng.
\newblock Training gans with optimism.
\newblock \emph{arXiv preprint arXiv:1711.00141}, 2017.

\bibitem[Dean et~al.(2012)Dean, Corrado, Monga, Chen, Devin, Mao, Ranzato, Senior, Tucker, Yang, et~al.]{dean2012large}
Jeffrey Dean, Greg Corrado, Rajat Monga, Kai Chen, Matthieu Devin, Mark Mao, Marc'aurelio Ranzato, Andrew Senior, Paul Tucker, Ke~Yang, et~al.
\newblock Large scale distributed deep networks.
\newblock \emph{Advances in Neural Information Processing Systems}, 25, 2012.

\bibitem[Defazio et~al.(2014{\natexlab{a}})Defazio, Bach, and Lacoste-Julien]{defazio2014saga}
Aaron Defazio, Francis Bach, and Simon Lacoste-Julien.
\newblock Saga: A fast incremental gradient method with support for non-strongly convex composite objectives.
\newblock \emph{Advances in Neural Information Processing Systems}, 27, 2014{\natexlab{a}}.

\bibitem[Defazio et~al.(2014{\natexlab{b}})Defazio, Domke, et~al.]{defazio2014finito}
Aaron Defazio, Justin Domke, et~al.
\newblock Finito: A faster, permutable incremental gradient method for big data problems.
\newblock In \emph{International Conference on Machine Learning}, pages 1125--1133. PMLR, 2014{\natexlab{b}}.

\bibitem[Esser et~al.(2010)Esser, Zhang, and Chan]{esser2010general}
Ernie Esser, Xiaoqun Zhang, and Tony~F Chan.
\newblock A general framework for a class of first order primal-dual algorithms for convex optimization in imaging science.
\newblock \emph{SIAM Journal on Imaging Sciences}, 3\penalty0 (4):\penalty0 1015--1046, 2010.

\bibitem[Facchinei and Pang(2003)]{facchinei2003finite}
Francisco Facchinei and Jong-Shi Pang.
\newblock \emph{Finite-dimensional variational inequalities and complementarity problems}.
\newblock Springer, 2003.

\bibitem[Fang et~al.(2018)Fang, Li, Lin, and Zhang]{fang2018spider}
Cong Fang, Chris~Junchi Li, Zhouchen Lin, and Tong Zhang.
\newblock Spider: Near-optimal non-convex optimization via stochastic path-integrated differential estimator.
\newblock \emph{Advances in Neural Information Processing Systems}, 31, 2018.

\bibitem[Gidel et~al.(2018)Gidel, Berard, Vignoud, Vincent, and Lacoste-Julien]{gidel2018variational}
Gauthier Gidel, Hugo Berard, Ga{\"e}tan Vignoud, Pascal Vincent, and Simon Lacoste-Julien.
\newblock A variational inequality perspective on generative adversarial networks.
\newblock \emph{arXiv preprint arXiv:1802.10551}, 2018.

\bibitem[Goodfellow et~al.(2014)Goodfellow, Pouget-Abadie, Mirza, Xu, Warde-Farley, Ozair, Courville, and Bengio]{goodfellow2014generative}
Ian Goodfellow, Jean Pouget-Abadie, Mehdi Mirza, Bing Xu, David Warde-Farley, Sherjil Ozair, Aaron Courville, and Yoshua Bengio.
\newblock Generative adversarial nets.
\newblock \emph{Advances in Neural Information Processing Systems}, 27, 2014.

\bibitem[Gorbunov et~al.(2020)Gorbunov, Hanzely, and Richt{\'a}rik]{gorbunov2020unified}
Eduard Gorbunov, Filip Hanzely, and Peter Richt{\'a}rik.
\newblock A unified theory of sgd: Variance reduction, sampling, quantization and coordinate descent.
\newblock In \emph{International Conference on Artificial Intelligence and Statistics}, pages 680--690. PMLR, 2020.

\bibitem[Gower et~al.(2020)Gower, Schmidt, Bach, and Richt{\'a}rik]{gower2020variance}
Robert~M Gower, Mark Schmidt, Francis Bach, and Peter Richt{\'a}rik.
\newblock Variance-reduced methods for machine learning.
\newblock \emph{Proceedings of the IEEE}, 108\penalty0 (11):\penalty0 1968--1983, 2020.

\bibitem[Gurbuzbalaban et~al.(2017)Gurbuzbalaban, Ozdaglar, and Parrilo]{gurbuzbalaban2017convergence}
Mert Gurbuzbalaban, Asuman Ozdaglar, and Pablo~A Parrilo.
\newblock On the convergence rate of incremental aggregated gradient algorithms.
\newblock \emph{SIAM Journal on Optimization}, 27\penalty0 (2):\penalty0 1035--1048, 2017.

\bibitem[G{\"u}rb{\"u}zbalaban et~al.(2021)G{\"u}rb{\"u}zbalaban, Ozdaglar, and Parrilo]{gurbuzbalaban2021random}
Mert G{\"u}rb{\"u}zbalaban, Asu Ozdaglar, and Pablo~A Parrilo.
\newblock Why random reshuffling beats stochastic gradient descent.
\newblock \emph{Mathematical Programming}, 186:\penalty0 49--84, 2021.

\bibitem[Haochen and Sra(2019)]{haochen2019random}
Jeff Haochen and Suvrit Sra.
\newblock Random shuffling beats sgd after finite epochs.
\newblock In \emph{International Conference on Machine Learning}, pages 2624--2633. PMLR, 2019.

\bibitem[Harker and Pang(1990)]{harker1990finite}
Patrick~T Harker and Jong-Shi Pang.
\newblock Finite-dimensional variational inequality and nonlinear complementarity problems: a survey of theory, algorithms and applications.
\newblock \emph{Mathematical Programming}, 48\penalty0 (1):\penalty0 161--220, 1990.

\bibitem[Hu et~al.(2019)Hu, Li, Lian, Liu, and Yuan]{hu2019efficient}
Wenqing Hu, Chris~Junchi Li, Xiangru Lian, Ji~Liu, and Huizhuo Yuan.
\newblock Efficient smooth non-convex stochastic compositional optimization via stochastic recursive gradient descent.
\newblock \emph{Advances in Neural Information Processing Systems}, 32, 2019.

\bibitem[Huang et~al.(2021)Huang, Yuan, Mao, and Yin]{huang2021improved}
Xinmeng Huang, Kun Yuan, Xianghui Mao, and Wotao Yin.
\newblock An improved analysis and rates for variance reduction under without-replacement sampling orders.
\newblock \emph{Advances in Neural Information Processing Systems}, 34:\penalty0 3232--3243, 2021.

\bibitem[Jin and Sidford(2020)]{jin2020efficiently}
Yujia Jin and Aaron Sidford.
\newblock Efficiently solving mdps with stochastic mirror descent.
\newblock In \emph{International Conference on Machine Learning}, pages 4890--4900. PMLR, 2020.

\bibitem[Johnson and Zhang(2013)]{johnson2013accelerating}
Rie Johnson and Tong Zhang.
\newblock Accelerating stochastic gradient descent using predictive variance reduction.
\newblock \emph{Advances in Neural Information Processing Systems}, 26, 2013.

\bibitem[Juditsky et~al.(2011)Juditsky, Nemirovski, and Tauvel]{juditsky2011solving}
Anatoli Juditsky, Arkadi Nemirovski, and Claire Tauvel.
\newblock Solving variational inequalities with stochastic mirror-prox algorithm.
\newblock \emph{Stochastic Systems}, 1\penalty0 (1):\penalty0 17--58, 2011.

\bibitem[Kinderlehrer and Stampacchia(2000)]{kinderlehrer2000introduction}
David Kinderlehrer and Guido Stampacchia.
\newblock \emph{An introduction to variational inequalities and their applications}.
\newblock SIAM, 2000.

\bibitem[Koloskova et~al.(2023)Koloskova, Doikov, Stich, and Jaggi]{koloskova2023convergence}
Anastasia Koloskova, Nikita Doikov, Sebastian~U Stich, and Martin Jaggi.
\newblock On convergence of incremental gradient for non-convex smooth functions.
\newblock \emph{arXiv preprint arXiv:2305.19259}, 2023.

\bibitem[Korpelevich(1976)]{korpelevich1976extragradient}
Galina~M Korpelevich.
\newblock The extragradient method for finding saddle points and other problems.
\newblock \emph{Matecon}, 12:\penalty0 747--756, 1976.

\bibitem[Kovalev et~al.(2020)Kovalev, Horv{\'a}th, and Richt{\'a}rik]{kovalev2020don}
Dmitry Kovalev, Samuel Horv{\'a}th, and Peter Richt{\'a}rik.
\newblock Don’t jump through hoops and remove those loops: Svrg and katyusha are better without the outer loop.
\newblock In \emph{Algorithmic Learning Theory}, pages 451--467. PMLR, 2020.

\bibitem[Li et~al.(2019)Li, Zhu, So, and Lee]{li2019incremental}
Xiao Li, Zhihui Zhu, Anthony Man-Cho So, and Jason~D Lee.
\newblock Incremental methods for weakly convex optimization.
\newblock \emph{arXiv preprint arXiv:1907.11687}, 2019.

\bibitem[Liang and Stokes(2019)]{liang2019interaction}
Tengyuan Liang and James Stokes.
\newblock Interaction matters: A note on non-asymptotic local convergence of generative adversarial networks.
\newblock In \emph{The 22nd International Conference on Artificial Intelligence and Statistics}, pages 907--915. PMLR, 2019.

\bibitem[Lu et~al.(2022)Lu, Meng, and De~Sa]{lu2022general}
Yucheng Lu, Si~Yi Meng, and Christopher De~Sa.
\newblock A general analysis of example-selection for stochastic gradient descent.
\newblock In \emph{International Conference on Learning Representations (ICLR)}, volume~10, 2022.

\bibitem[Madry et~al.(2017)Madry, Makelov, Schmidt, Tsipras, and Vladu]{madry2017towards}
Aleksander Madry, Aleksandar Makelov, Ludwig Schmidt, Dimitris Tsipras, and Adrian Vladu.
\newblock Towards deep learning models resistant to adversarial attacks.
\newblock \emph{arXiv preprint arXiv:1706.06083}, 2017.

\bibitem[Malinovsky et~al.(2023)Malinovsky, Sailanbayev, and Richt{\'a}rik]{malinovsky2023random}
Grigory Malinovsky, Alibek Sailanbayev, and Peter Richt{\'a}rik.
\newblock Random reshuffling with variance reduction: New analysis and better rates.
\newblock In \emph{Uncertainty in Artificial Intelligence}, pages 1347--1357. PMLR, 2023.

\bibitem[Malitsky(2015)]{malitsky2015projected}
Yu~Malitsky.
\newblock Projected reflected gradient methods for monotone variational inequalities.
\newblock \emph{SIAM Journal on Optimization}, 25\penalty0 (1):\penalty0 502--520, 2015.

\bibitem[Malitsky and Tam(2020)]{malitsky2020forward}
Yura Malitsky and Matthew~K Tam.
\newblock A forward-backward splitting method for monotone inclusions without cocoercivity.
\newblock \emph{SIAM Journal on Optimization}, 30\penalty0 (2):\penalty0 1451--1472, 2020.

\bibitem[Mangasarian and Solodov(1993)]{mangasarian1993serial}
Olvi~L Mangasarian and MV~Solodov.
\newblock Serial and parallel backpropagation convergence via nonmonotone perturbed minimization.
\newblock Technical report, University of Wisconsin-Madison Department of Computer Sciences, 1993.

\bibitem[Medyakov et~al.(2023)Medyakov, Molodtsov, Beznosikov, and Gasnikov]{medyakov2023optimal}
Daniil Medyakov, Gleb Molodtsov, Aleksandr Beznosikov, and Alexander Gasnikov.
\newblock Optimal data splitting in distributed optimization for machine learning.
\newblock In \emph{Doklady Mathematics}, volume 108, pages S465--S475. Springer, 2023.

\bibitem[Mertikopoulos et~al.(2018)Mertikopoulos, Lecouat, Zenati, Foo, Chandrasekhar, and Piliouras]{mertikopoulos2018optimistic}
Panayotis Mertikopoulos, Bruno Lecouat, Houssam Zenati, Chuan-Sheng Foo, Vijay Chandrasekhar, and Georgios Piliouras.
\newblock Optimistic mirror descent in saddle-point problems: Going the extra (gradient) mile.
\newblock \emph{arXiv preprint arXiv:1807.02629}, 2018.

\bibitem[Mishchenko et~al.(2020{\natexlab{a}})Mishchenko, Khaled, and Richt{\'a}rik]{mishchenko2020random}
Konstantin Mishchenko, Ahmed Khaled, and Peter Richt{\'a}rik.
\newblock Random reshuffling: Simple analysis with vast improvements.
\newblock \emph{Advances in Neural Information Processing Systems}, 33:\penalty0 17309--17320, 2020{\natexlab{a}}.

\bibitem[Mishchenko et~al.(2020{\natexlab{b}})Mishchenko, Kovalev, Shulgin, Richt{\'a}rik, and Malitsky]{mishchenko2020revisiting}
Konstantin Mishchenko, Dmitry Kovalev, Egor Shulgin, Peter Richt{\'a}rik, and Yura Malitsky.
\newblock Revisiting stochastic extragradient.
\newblock In \emph{International Conference on Artificial Intelligence and Statistics}, pages 4573--4582. PMLR, 2020{\natexlab{b}}.

\bibitem[Mohtashami et~al.(2022)Mohtashami, Stich, and Jaggi]{mohtashami2022characterizing}
Amirkeivan Mohtashami, Sebastian Stich, and Martin Jaggi.
\newblock Characterizing \& finding good data orderings for fast convergence of sequential gradient methods.
\newblock \emph{arXiv preprint arXiv:2202.01838}, 2022.

\bibitem[Mokhtari et~al.(2018)Mokhtari, Gurbuzbalaban, and Ribeiro]{mokhtari2018surpassing}
Aryan Mokhtari, Mert Gurbuzbalaban, and Alejandro Ribeiro.
\newblock Surpassing gradient descent provably: A cyclic incremental method with linear convergence rate.
\newblock \emph{SIAM Journal on Optimization}, 28\penalty0 (2):\penalty0 1420--1447, 2018.

\bibitem[Mokhtari et~al.(2020)Mokhtari, Ozdaglar, and Pattathil]{mokhtari2020unified}
Aryan Mokhtari, Asuman Ozdaglar, and Sarath Pattathil.
\newblock A unified analysis of extra-gradient and optimistic gradient methods for saddle point problems: Proximal point approach.
\newblock In \emph{International Conference on Artificial Intelligence and Statistics}, pages 1497--1507. PMLR, 2020.

\bibitem[Moulines and Bach(2011)]{moulines2011non}
Eric Moulines and Francis Bach.
\newblock Non-asymptotic analysis of stochastic approximation algorithms for machine learning.
\newblock \emph{Advances in Neural Information Processing Systems}, 24, 2011.

\bibitem[Nagaraj et~al.(2019)Nagaraj, Jain, and Netrapalli]{nagaraj2019sgd}
Dheeraj Nagaraj, Prateek Jain, and Praneeth Netrapalli.
\newblock Sgd without replacement: Sharper rates for general smooth convex functions.
\newblock In \emph{International Conference on Machine Learning}, pages 4703--4711. PMLR, 2019.

\bibitem[Nedic and Bertsekas(2001)]{nedic2001incremental}
Angelia Nedic and Dimitri~P Bertsekas.
\newblock Incremental subgradient methods for nondifferentiable optimization.
\newblock \emph{SIAM Journal on Optimization}, 12\penalty0 (1):\penalty0 109--138, 2001.

\bibitem[Nemirovski(2004)]{nemirovski2004prox}
Arkadi Nemirovski.
\newblock Prox-method with rate of convergence o (1/t) for variational inequalities with lipschitz continuous monotone operators and smooth convex-concave saddle point problems.
\newblock \emph{SIAM Journal on Optimization}, 15\penalty0 (1):\penalty0 229--251, 2004.

\bibitem[Nesterov(2005)]{nesterov2005smooth}
Yu~Nesterov.
\newblock Smooth minimization of non-smooth functions.
\newblock \emph{Mathematical Programming}, 103:\penalty0 127--152, 2005.

\bibitem[Nesterov(2007)]{nesterov2007dual}
Yurii Nesterov.
\newblock Dual extrapolation and its applications to solving variational inequalities and related problems.
\newblock \emph{Mathematical Programming}, 109\penalty0 (2):\penalty0 319--344, 2007.

\bibitem[Nguyen et~al.(2017)Nguyen, Liu, Scheinberg, and Tak{\'a}{\v{c}}]{nguyen2017sarah}
Lam~M Nguyen, Jie Liu, Katya Scheinberg, and Martin Tak{\'a}{\v{c}}.
\newblock Sarah: A novel method for machine learning problems using stochastic recursive gradient.
\newblock In \emph{International Conference on Machine Learning}, pages 2613--2621. PMLR, 2017.

\bibitem[Omidshafiei et~al.(2017)Omidshafiei, Pazis, Amato, How, and Vian]{omidshafiei2017deep}
Shayegan Omidshafiei, Jason Pazis, Christopher Amato, Jonathan~P How, and John Vian.
\newblock Deep decentralized multi-task multi-agent reinforcement learning under partial observability.
\newblock In \emph{International Conference on Machine Learning}, pages 2681--2690. PMLR, 2017.

\bibitem[Palaniappan and Bach(2016)]{palaniappan2016stochastic}
Balamurugan Palaniappan and Francis Bach.
\newblock Stochastic variance reduction methods for saddle-point problems.
\newblock \emph{Advances in Neural Information Processing Systems}, 29, 2016.

\bibitem[Peng et~al.(2020)Peng, Dai, Zhang, and Cheng]{peng2020training}
Wei Peng, Yu-Hong Dai, Hui Zhang, and Lizhi Cheng.
\newblock Training gans with centripetal acceleration.
\newblock \emph{Optimization Methods and Software}, 35\penalty0 (5):\penalty0 955--973, 2020.

\bibitem[Rajput et~al.(2020)Rajput, Gupta, and Papailiopoulos]{rajput2020closing}
Shashank Rajput, Anant Gupta, and Dimitris Papailiopoulos.
\newblock Closing the convergence gap of sgd without replacement.
\newblock In \emph{International Conference on Machine Learning}, pages 7964--7973. PMLR, 2020.

\bibitem[Robbins and Monro(1951)]{robbins1951stochastic}
Herbert Robbins and Sutton Monro.
\newblock A stochastic approximation method.
\newblock \emph{The Annals of Mathematical Statistics}, pages 400--407, 1951.

\bibitem[Roux et~al.(2012)Roux, Schmidt, and Bach]{roux2012stochastic}
Nicolas Roux, Mark Schmidt, and Francis Bach.
\newblock A stochastic gradient method with an exponential convergence rate for finite training sets.
\newblock \emph{Advances in Neural Information Processing Systems}, 25, 2012.

\bibitem[Safran and Shamir(2020)]{safran2020good}
Itay Safran and Ohad Shamir.
\newblock How good is sgd with random shuffling?
\newblock In \emph{Conference on Learning Theory}, pages 3250--3284. PMLR, 2020.

\bibitem[Stich(2019)]{stich2019unified}
Sebastian~U Stich.
\newblock Unified optimal analysis of the (stochastic) gradient method.
\newblock \emph{arXiv preprint arXiv:1907.04232}, 2019.

\bibitem[Tseng(2000)]{tseng2000modified}
Paul Tseng.
\newblock A modified forward-backward splitting method for maximal monotone mappings.
\newblock \emph{SIAM Journal on Control and Optimization}, 38\penalty0 (2):\penalty0 431--446, 2000.

\bibitem[Von~Neumann and Morgenstern(1953)]{von1953theory}
John Von~Neumann and Oskar Morgenstern.
\newblock \emph{Theory of games and economic behavior: by J. Von Neumann and O. Morgenstern}.
\newblock Princeton University Press, 1953.

\bibitem[Xu et~al.(2004)Xu, Neufeld, Larson, and Schuurmans]{xu2004maximum}
Linli Xu, James Neufeld, Bryce Larson, and Dale Schuurmans.
\newblock Maximum margin clustering.
\newblock \emph{Advances in Neural Information Processing Systems}, 17, 2004.

\bibitem[Ying et~al.(2020)Ying, Yuan, and Sayed]{ying2020variance}
Bicheng Ying, Kun Yuan, and Ali~H Sayed.
\newblock Variance-reduced stochastic learning under random reshuffling.
\newblock \emph{IEEE Transactions on Signal Processing}, 68:\penalty0 1390--1408, 2020.

\end{thebibliography}
\end{document}